\documentclass[a4paper,11pt]{article}
\usepackage[top=25mm,left=22mm]{geometry}
\usepackage{amsmath,amssymb,amsthm,amsfonts}
\usepackage{lineno,url}\usepackage{float}
\floatstyle{plaintop}\restylefloat{table}
\usepackage[tableposition=top]{caption}
\usepackage{graphicx}
\usepackage{multirow}
\usepackage{pstricks}
\usepackage{fancyvrb}
\usepackage{paralist}
\newtheorem{theorem}{Theorem}[section]
\newtheorem{remark}[theorem]{Remark}
\def\bexe{\begin{exercise}}\def\eexe{\eex\end{exercise}}
\def\bsol{\begin{solution}}\def\esol{\eex\end{solution}}
\def\bexa{\begin{example}}\def\eexa{\end{example}}
\def\brem{\begin{remark}}\def\erem{\end{remark}}
\def\bthm{\begin{theorem}}\def\ethm{\end{theorem}}
\def\blem{\begin{lemma}}\def\elem{\end{lemma}}
\def\bcor{\begin{corollary}}\def\ecor{\end{corollary}}
\def\bdefi{\begin{definition}}\def\edefi{\end{definition}}
\newcommand{\IDEA}{\textbf{Idea of the Proof.} }

\def\bmip{\begin{minipage}{\textwidth}}\def\emip{\end{minipage}}
\def\huga#1{\begin{gather} #1 \end{gather}}

\def\hual#1{\begin{align} #1 \end{align}}

\newcommand{\R}{{\mathbb R}}
\newcommand{\C}{{\mathbb C}}

\def\CD{{\cal D}}  
\def\CH{{\cal H}}

\def\uti{\tilde{u}}

\def\ga{\gamma}\def\uh{\hat{u}}\def\vh{\hat{v}}
\def\noi{\noindent}\def\ds{\displaystyle}
\def\pa{{\partial}}\def\lam{\lambda}
\newcommand{\bi}{\begin{itemize}}\newcommand{\ei}{\end{itemize}}
\newcommand{\ben}{\begin{enumerate}}\newcommand{\een}{\end{enumerate}}
\newcommand{\bce}{\begin{center}}\newcommand{\ece}{\end{center}}
\newcommand{\bci}{\begin{compactitem}}\newcommand{\eci}{\end{compactitem}}
\newcommand{\bcen}{\begin{compactenum}}\newcommand{\ecen}{\end{compactenum}}
\newcommand{\reff}[1]{(\ref{#1})}
\newcommand{\ov}[1]{{\overline {#1}}}

\newcommand{\vs}[1]{{\vspace{#1}}}

\def\ra{\rightarrow}

\newcommand{\barr}{\begin{array}}\newcommand{\earr}{\end{array}}
\newcommand{\bpm}{\begin{pmatrix}}\newcommand{\epm}{\end{pmatrix}}
\newcommand{\bsm}{\left(\begin{smallmatrix}}
\newcommand{\esm}{\end{smallmatrix}\right)}
\newcommand{\ba}{\begin{array}}\newcommand{\ea}{\end{array}}
\def\dd{\, {\rm d}}

\def\er{{\rm e}}

\def\Om{\Omega}
\def\ddt{\frac{\rm d}{{\rm d}t}}
\def\del{\delta}

\def\eex{\hfill\mbox{$\rfloor$}}

\def\sig{\sigma}
\def\al{\alpha}

\def\Lam{\Lambda}

\def\bd{\begin{displaymath}} \def\ed{\end{displaymath}}
\def\ba{\begin{array}} \def\ea{\end{array}}

\newcommand{\argmax}{\operatornamewithlimits{argmax}}
\def\PMAXP{Pontryagin's Maximum Principle}

\def\ig{\includegraphics}\def\pdep{{\tt pde2path}}\def\mtom{{\tt mtom}}
\def\poc{{\tt p2poc}}
\def\pocl{{\tt p2poclib}}
\def\mlab{{\tt Matlab}}\def\ptool{{\tt pdetoolbox}}

\def\Dx{\dd x}\def\evalat#1{#1}
\renewcommand{\arraystretch}{1}\renewcommand{\baselinestretch}{1.0}
\def\medskip{}\def\bigskip{}

\usepackage{hyperref}

\begin{document}
\mbox{}\vspace{0.1cm}\begin{center}\Large
The \pdep\ add-on toolbox \poc\ for solving infinite time--horizon 
spatially distributed optimal control problems\\[2mm]
\normalsize -- Quickstart Guide --\\[2mm] 
 Hannes Uecker$^1$\\[2mm]
\footnotesize
$^1$ Institut f\"ur Mathematik, Universit\"at Oldenburg, D26111 Oldenburg, 
hannes.uecker@uni-oldenburg.de\\[2mm]
\normalsize
\today
\end{center}
\noi
\begin{abstract} 
\poc\ is an add--on toolbox to the \mlab\ package \pdep. It is aimed at the 
numerical solution of infinite time horizon 
optimal control (OC) problems for 
parabolic systems of PDE over 1D or 2D spatial domains. The basic idea is to 
treat the OC problem via the associated canonical system in two steps. 
First we use \pdep\ to find branches of stationary solutions 
of the canonical system, also called canonical steady states (CSS). 
In a second step 
we use the results and the spatial discretization of the first step 
to calculate the objective values of time-dependent canonical paths 
ending at a CSS with the so called saddle point property. 
This is a (typically very high 
dimensional) boundary value problem (BVP) in time, which we solve 
by combining a modification of the BVP solver TOM with a continuation 
algorithm in the initial states. 
We explain the design and usage of the package via two example problems, 
namely the optimal management of a distributed shallow lake 
model, and of a semi-arid grazing system. Both show interesting 
bifurcations of so called patterned CSS, and in particular 
the latter also a variety 
of patterned {\em optimal} steady states. 
The package (library and demos) can be downloaded at 
{\tt www.staff.uni-oldenburg.de/hannes.uecker/pde2path}. 
\end{abstract}
\noindent


\section{Introduction}\label{i-sec}
Denoting the state variable (vector) by $v=v(t,x)\in\R^N$, and the 
control by $k=k(t,x)\in\R$,   
we consider spatially distribute infinite time horizon 
optimal control (OC) problems of the form 
\begin{subequations} \label{oc1}
  \begin{align}
    &V(v_0(\cdot)):=\max_{k(\cdot,\cdot)}J(v_0(\cdot),k(\cdot,\cdot)), \qquad 
    J(v_0(\cdot),k(\cdot,\cdot)):=\int_0^\infty\er^{-\rho t}
J_{ca}(v(t),k(t)) \dd t,\\
\text{where }& J_{ca}(v(\cdot,t),k(\cdot,t))=\frac 1{|\Om|}
\int_\Om J_c(v(x,t),k(x,t))\dd x\\
\intertext{is the spatially averaged current value objective function, 
with the local current value $J_c{:}\R^{N+1}{\ra}\R$ a given function; 
$\rho>0$ is the discount rate, 
and $v:\Om\times [0,\infty)\ra\R^N$ fulfills a PDE of the form} 
&\pa_t v=-G_1(v,k):=D\Delta v+g_1(v,k), \quad v|_{t=0}=v_0. 
  \end{align}
\end{subequations}
Here, $D\in\R^{N\times N}$ is a diffusion matrix, 
$\Delta=\pa_{x_1}^2+\ldots+\pa_{x_d}^2$ is the Laplacian, and 
(\ref{oc1}c) holds in a bounded domain $\Om\subset\R^d$ with suitable boundary 
conditions (BC), where for simplicity we restrict to homogeneous Neumann 
$\pa_\nu v=0$, $\nu$ the outer normal.  
In applications, $J_c$ and $G_1$ of course often also 
depend on a number of parameters, 
which however for simplicity we do not display 
here.\footnote{$G_1$ in (\ref{oc1}c) can in fact be of a much more general form, 
but for simplicity here we stick to  (\ref{oc1}c). The convention 
that $\pa_t v=-G_1(v)$ instead of $\pa_t v=G_1(v)$ with 
$G_1(v)=D\Delta v+g_1(v)$ is inherited via \pdep\ from the \mlab\ \ptool, which 
assembles $-\Delta$ into the stiffness matrix $K$.}

Introducing the costates $\lam:\Om\times(0,\infty)\ra \R^{N}$ and the 
(local current value) Hamiltonian
\hual{
\CH&=\CH(v,\lam,k)=J_c(v,k)+\lam^T(D\Delta v+g_1(v,k)),} 
by Pontryagin's Maximum Principle (see the references below) 
for $\tilde{\CH}=
\int_0^\infty \er^{-\rho t} \ov{\CH}(t)\dd t$ with the spatial integral 
\huga{\label{fullH} 
\ov{\CH}(t)=\int_\Om \CH(v(x,t),\lam(x,t),k(x,t))
\dd x, 
}
an optimal solution $(v,\lam)$ (or equivalently 
$(v,k):\Om\times[0,\infty)\ra \R^{N+1})$) has to solve the canonical system (CS) 
\begin{subequations} \label{cs}
\hual{
\pa_t v&=\pa_\lam\CH=D\Delta v+g_1(v,k), \quad v|_{t=0}=v_0, \\
\pa_t \lam&=\rho\lam-\pa_v\CH=\rho\lam+g_2(v,k)-D\Delta\lam, 
\intertext{
where  $k=\argmax_{\tilde{k}}\CH(v,\lam,\tilde{k})$, which generally we assume to 
be obtained from solving}
\pa_k\CH(v,\lam,k)&=0.
}
\end{subequations}
The costates $\lam$ also fulfill zero flux BC, and 
derivatives like $\pa_v \CH$ etc are taken 
variationally, i.e., for $\ov{\CH}$. For instance, for $N=1$ and 
$\Phi(v,\lam):=\lam\Delta v$ we have $\ov{\Phi}(v,\lam)
=\int_\Om \lam\Delta v\dd x
=\int_\Om (\Delta\lam)v\dd x$ by Gau\ss' theorem, hence 
$\delta_v \ov{\Phi}(v,\lam)[h]=\int (\Delta\lam) h\dd x$, and 
by the Riesz representation theorem we 
identify $\delta_v \ov{\Phi}(v,\lam)$ and hence $\pa_v\Phi(v,\lam)$ 
with the multiplier $\Delta\lam$.
(\reff{cs}c) typically applies and yields a unique solution $k$ under suitable 
concavity assumptions on $\CH$, and in the absence of control constraints, 
see below. 

In principle we want so solve \reff{cs} for $t\in[0,\infty)$, but 
in (\ref{cs}a) we have initial data for only half the variables, 
and in (\ref{cs}b) we have anti--diffusion, such that \reff{cs} 
is ill--posed as an initial value problem. 
For convenience we set 
\huga{
u(t,\cdot):=\bpm v(t,\cdot)\\ \lam(t,\cdot)\epm: \Om\ra\R^{2N}, 
} 
and write \reff{cs} as 
\begin{subequations} \label{cs2}
 \hual{\label{cs2a}
\pa_t u&=-G(u,\eta):=\CD\Delta u+f(u), \text{ where } 
\CD=\bpm D&0\\0&-D\epm, \quad f(u)=\bpm g_1(u)\\ g_2(u)\epm, 
\intertext{and where $\eta\in\R^p$ stands for parameters present, which for 
instance also includes the discount rate $\rho$. Besides the 
boundary condition $\pa_\nu u=0$ and the initial 
condition }
v|_{t=0}&=v_0, 
\intertext{we then impose 
the tranversality condition} 
\lim_{t\to\infty}\er^{-\rho t}u(t)&=0. \label{tcond}
}
\end{subequations}
A solution $u$ of the canonical system \reff{cs2} is called a 
\emph{canonical path}, 
and an equilibrium of \reff{cs2a} (which automatically fulfils \reff{tcond}) 
is called 
\emph{canonical steady state (CSS)}. With a slight abuse of notation 
we also call $(v,k)$ with $k$ given by (\ref{cs}d) a canonical 
path.  See also, e.g., \cite{GU15} for more formal definitions, 
further comments on the notions of optimal systems, and, e.g., the significance 
of the transversality condition \reff{tcond}.

For general background on OC in a PDE setting see \cite{Tr10} and the 
references therein, or specifically \cite{RZ99, RZ99b} and 
\cite[Chapter5]{AAC11} for \PMAXP{} 
for OC problems for semilinear diffusive models. However, these works are in a 
finite time horizon setting, and often the objective function 
is linear in the control and there are control constraints, e.g., 
$k(x,t)\in K$ with some bounded interval $K$. Therefore $k$ is not obtained from 
the analogue of (\ref{cs}c), but rather takes the values from $\pa K$, which 
is often called bang control.  
In, e.g., \cite{CP12, ACKT13}, some specific models have been studied 
in this setting and a rather theoretical way, i.e., the focus is on 
deriving the canonical system and showing well-posedness and 
the existence of an optimal control. \cite{Apre14} additionally 
contains numerical simulations for a finite time horizon control--constrained 
OC problem for a three species spatial predator-prey system, again 
leading to bang type controls.  See also \cite{NPS11} 
and the references therein for numerical methods for (finite time horizon)  
constrained parabolic optimal control problems.

Here we do not (yet) consider (active) 
control or state constraints, and no terminal time, but the infinite 
time horizon. 
Our models and method are motivated by \cite{BX08, BX10}, which also 
discuss \PMAXP{} in this setting.  We do not extend the theory, but 
rather consider 
\reff{cs2} after a spatial discretization as a (large) ODE problem, 
and essentially treat this using the notations and ideas from 
\cite{grassetal2008}, to give a numerical framework to calculate 
optimal solutions. 
Using the canonical system \reff{cs} 
we proceed in two steps, which can be seen as a variant of the ``connecting 
orbit method'', see, e.g., \cite{BPS01}, and 
also \S\ref{dsec} for further background and remarks on 
the related literature: first we compute (branches of) CSS, and 
second we compute canonical paths connecting to some CSS.  
This also means that we take a somewhat 
broader perspective than aiming at computing just one optimal control, 
given an initial condition $v_0$, which without further information 
is ill-posed anyway. Instead, our method aims to give a somewhat global 
picture by identifying the pertinent CSS and their respective domains 
of attraction. 
\paragraph{(a) CSS branches.} 
We compute (approximate) CSS of \reff{cs2}, i.e., solutions $\uh$ of 
\huga{\label{cs21}
G(u,\eta)=0, 
}
together with the BC. 
For this use the package \pdep\ \cite{p2p, p2p2} to set up a FEM 
discretization of \reff{cs21} as a continuation/bifurcation problem 
in one of the parameters, which we call $\eta$ again. This gives 
branches $\eta\mapsto \uh(\eta)$ of solutions, which is in particular 
useful to possibly find 
several solutions $\uh^{(l)}(\eta)$, $j=l,\ldots, m$ at fixed $\eta$. 
By computing the associated $J_{ca}(\hat v^{(l)},k^{(l)})$ 
we can identify which of these is optimal amongst the CSS. 
Given a CSS $\uh$, for simplicity we also write 
$J_{ca}(\uh):=J_{ca}(\hat v^{(l)},k^{(l)})$, and moreover, have 
\huga{J(\uh)=J_{ca}(\uh)/\rho.
}

\paragraph{(b) Canonical paths.}
In a second step (b),  we calculate canonical 
paths ending at a CSS $\uh$ (and 
often starting at the state values $\vh_0$ of a different CSS $\uh_0$), 
and the objective 
values of the canonical paths. For this we choose a truncation time 
$T$ and modify \reff{tcond} to the condition 
that $u(T)\in W_s(\uh)$ and near $\uh$, where $W_s(\uh)$ denotes the stable 
manifold of $\uh$. In practice, we approximate $W_s(\uh)$ by the stable 
eigenspace $E_s(\uh)$, and thus consider the BVP 
\begin{subequations}\label{bvp0}
\hual{
\pa_t u&=-G(u),\\
v|_{t=0}&=v_0,\\
u(T)&\in E_s(\uh) \text{ (and $\|u(T)-\uh\|$ small)}. 
} 
\end{subequations}
Using the spatial FEM discretizations, the implementation 
of $G$, and the results from the first step (a), 
if the mesh in the FEM consists of $n$ nodes, then 
$u(t)\in\R^{2Nn}$, and (\ref{bvp0}a) yields a system of $2Nn$ ODEs 
in the form (with a slight abuse of notation) 
\begin{subequations}\label{bvp1}
\hual{
M\ddt u&=-G(u),
\intertext{while the initial and transversality conditions become}
v|_{t=0}&=v_0,\\
\Psi(u(T)-\uh)&=0 \text{ (and  $\|u(T)-\uh\|$ small)}. 
} 
\end{subequations}
Here $M\in\R^{2Nn\times 2Nn}$ is the mass matrix of the FEM mesh, and 
$\Psi\in \R^{Nn\times 2Nn}$ defines the projection onto $E_u(\uh)$. Moreover, 
(\ref{bvp1}b) consists of $Nn$ initial conditions for 
the states, while the costates $\lam$ (and hence the control $k$) 
are free. Thus, to have $2Nn$ BC altogether we need dim$E_s(\uh)=Nn$. 
On the other hand, we always have dim$E_s(\uh)\le Nn$, see 
\cite[Appendix A]{GU15}. We define the defect 
\huga{\label{ddef}
d(\uh):={\rm dim}E_s(\uh)-Nn}
and call a CSS $\uh$ with $d(\uh)=0$ a CSS with the saddle--point 
property (SPP). At first sight it may appear that $d(\uh)$ depends 
on the spatial discretization, i.e., on the number of $n$ of nodes. 
However, $d(\uh)$ remains constant for finer and finer meshes, 
see  \cite[Appendix A]{GU15} for further comments. 

For $\uh=(\vh,\hat\lam)$ with the SPP, 
and $\|v_0-\vh\|$ sufficiently small, we may expect the existence 
of a solution $u$ of \reff{bvp1}, which moreover can be found from 
a Newton loop for \reff{bvp1} with initial guess $u(t)\equiv \uh$. 
On the other hand, for larger $\|v_0-\vh\|$ a solution of \reff{bvp1} 
may not exist, or a good initial guess may be hard 
to find, and therefore we use a continuation process also for \reff{bvp1}. 
In the simplest setting, assume that for some $\al\in[0,1]$
we have a solution $u_\al$ of \reff{bvp1} with (\ref{bvp1}b) replaced by 
\huga{\label{albc}
	v(0)=\al v_{0}+(1-\al)\hat v,
}
 (e.g., $\al=0$ and $u\equiv \uh$). We then increase $\al$ by some stepsize 
$\del_\al$ and use $u_\al$ as initial guess for (\ref{bvp1}a), 
(\ref{bvp1}c) and \reff{albc}, ultimately aiming at $\al=1$. 

To actually solve (\ref{bvp1}a), (\ref{bvp1}c) and \reff{albc} 
we use  TOM \cite{mazS2002,MT04,MT2009} (see also 
\url{www.dm.uniba.it/~mazzia/mazzia/?page_id=433}) 
in a version {\tt mtom} which accounts for the mass matrix $M$ on the 
lhs of (\ref{bvp1}a).\footnote{mtom is an ad-hoc modification of 
TOM, and  will be replaced by an official version of TOM which 
handles mass matrices once that becomes available. 
In the following, when discussing, 
e.g., the behaviour and usage of {\tt mtom}, we note that 
almost all of this is derived from TOM.}  
This  predictor ($u_\al$) -- corrector ({\tt mtom} for $\al+\del_\al$) 
continuation method corresponds to the ``natural'' parametrization 
of the continuation by $\al$, and is thus 
implemented in \poc\ as {\tt iscnat} (Initial State Continuation NATural). 
We also give the option to use a secant predictor 
 \huga{\label{undef0}
u^j(t)=u^{j-1}(t)+\del_\al \tau(t), \quad 
\tau(t)=\bigl(u^{(j-1)}(t)-u^{(j-2)}(t)\bigr)
/\|u^{(j-1)}(\cdot)-u^{(j-2)}(\cdot)\|_2, 
} 
where $u^{j-2}$ and $u^{j-1}$ are the two previous steps. However, 
the corrector still works at fixed $\al$, in contrast to 
the arclength predictor--corrector {\tt iscarc} described next. 

It may happen that no solution of 
(\ref{bvp1}a), (\ref{bvp1}c) and \reff{albc} 
 is found for $\al>\al_0$
for some $\al_0<1$, i.e., that the continuation to the intended
initial states fails. In that case, often the BVP shows a
fold in $\al$, and we use a modified continuation process, 
letting $\alpha$ be a free parameter and using a pseudo--arclength 
parametrization by $\sig$ in the BC at $t=0$. 
 Since {\tt mtom} does not allow free parameters we add the dummy ODE
$\dot\alpha=0$, and BCs at continuation step $j$,  
 \huga{\label{iccarc}
	\bigl\langle s,(u(0)-u^{(j-1)}(0)\rangle
+s_\alpha(\alpha-\alpha^{(j-1)})=\sig,
}
      with $u^{(j-1)}(\cdot)$ the solution from the
      previous continuation step $j-1$, and 
$(s,s_\alpha)\in\R^{2N}\times \R$ appropriately
      chosen with $\|(s,s_\alpha)\|_*=1$, where $\|\cdot\|_*$ is a 
suitable norm in $\R^{2N+1}$, which may contain different 
weights of $v$ and $v_\alpha$.  For $s=0$ and
      $s_\alpha=1$ we find {\tt iscnat} with stepsize $\del_\al=\sig$ again. 
To get around folds we may use the secant
$$s:=
\xi\bigl(u^{(j-1)}(0)-u^{(j-2)}(0)\bigr)/\|u^{(j-1)}(0)-u^{(j-2)}(0)\|_2
\text{ and } s_\al=1-\xi 
$$
with small $\xi$, and also a secant predictor 
\huga{\label{undef}
(u^j,\al^j)^{\text{pred}}
=(u^{j-1},\al^{j-1})+\sig \tau} 
for $t\mapsto u^j(t)$ with 
\huga{\label{spred} 
\tau=\xi\bigl(u^{(j-1)}(\cdot)-u^{(j-2)}(\cdot)\bigr)
/\|u^{(j-1)}(\cdot)-u^{(j-2)}(\cdot)\|_2
\text{ and } \tau_\al=1-\xi. 
}
This essentially follows \cite[\S7.2]{grassetal2008}, and is implemented 
in a routine {\tt iscarc}  (Initial State Continuation ARClength).

Finally, given $\uh$, to calculate $\Psi$, at startup we 
solve the generalized adjoint eigenvalue problem 
\huga{\label{psievp}
\pa_u G(\uh)^T\Phi=\Lambda M \Phi
} for the eigenvalues $\Lambda$ 
and (adjoint) eigenvectors $\Phi$, which also gives the defect $d(\uh)$ 
by counting the negative eigenvalues in $\Lam$. If 
$d(\uh)=0$, then from $\Phi\in \C^{2Nn\times 2Nn}$ 
we generate a real base 
of $E_u(\hat u)$ which we sort into the matrix  $\Psi\in \R^{Nn\times 2Nn}$. \\[2mm]

\noindent 
{\bf Acknowledgement.} I thank D. Grass, ORCOS Wien, for introducing 
me to the field of optimal control over infinite time horizons, and 
 for clarifying (and posing) many questions regarding the aim of software 
based on Pontryagin's Maximum Principle for spatially distributed OC problems.

\section{Examples and implementation details}
\subsection{The SLOC model}\label{slsec}
Following \cite{BX08}, in \cite{GU15} we consider a model 
for phosphorus $P=P(t,x)$ in a shallow lake, and phosphate load $k=k(x,t)$ 
as a control, which in 0D, i.e., in the ODE setting, has analyzed in detail 
for instance in \cite{KW10}. Here we explain how we set up the spatial 
so called 
Shallow Lake Optimal Control (SLOC) problem in \poc, and 
refer to \cite{GU15} for details about 
the modelling and the interpretation of results. The model reads 
\begin{subequations} 
\label{sldiff1}
  \begin{align}
    &V(P_0(\cdot)):=\max_{k(\cdot,\cdot)}J(P_0(\cdot),k(\cdot,\cdot)), \qquad 
    J(P_0(\cdot),k(\cdot,\cdot)):=\int_0^\infty\er^{-\rho t}
J_{ca}(P(t),k(t)), 
\dd t\label{sldiffusion_model_obj2}\\
\intertext{where $J_c(P,k)=\ln(k)-\ga P^2$ is the {\em local} current value 
objective function,}& J_{ca}(P(\cdot,t),k(\cdot,t))=\frac 1 {|\Om|}
\int_\Om J_c(P(x,t),k(x,t))\!\Dx\\
\intertext{is the spatially averaged current value objective function, and 
$P$ fulfills the PDE 
}
 &\pa_t P(x,t)=k(x,t)-bP(x,t)+\frac{P(x,t)^2}{1+P(x,t)^2}+D\Delta P(x,t),\label{sldiffusion_model_dyn}\\
    &\pa_\nu P(x,t)_{\pa\Omega}=0,\label{sldiffusion_model_bc}\qquad 
    P(x,t)_{t=0}=P_0(x),\quad
    x\in\Om\subset\R^d. 
  \end{align}
\end{subequations}
The parameter $b>0$ is the phosphor degradation rate, and $\ga>0$ 
are ecological costs of the phosphor contamination $P$. 
One wants a low $P$ for ecological reasons, but for 
economical reasons a high phosphate load 
$k$, for instance from fertilizers used by farmers. Thus, 
the objective function consists of the concave
increasing function $\ln(k)$, and the concave decreasing function
$-\ga P^2$.  We consider two scenarios, namely 
\hual{\text{ Scenario 1: } D=0.5, \rho=0.03, \ga=0.5, b\in (0.5,0.8) \text{ (primary
    bif.~param.)},\label{parsel}\\
\text{ Scenario 2: } D=0.5, \rho=0.3, b=0.55, \ga\in (2.5,3.7) \text{ (primary
    bif.~param.)}. 
}
With the co-state $q$ and local current value Hamiltonian 
\huga{\label{hpde}
\CH(P,q,\lam)=J_{c}(P,k)+q\bigl[k-bP+\frac{P^2}{1+P^2}+D\Delta P\bigr], 
} the canonical system for \reff{sldiff1} becomes, 
with $\ds k(x,t)=-\frac{1}{q(x,t)}$, 
\begin{subequations}
  \label{slcan2}
  \begin{align}
    {\pa_t} P(x,t)&=k(x,t)-bP(x,t)+\frac{P(x,t)^2}{1+P(x,t)^2}
+D\Delta P(x,t),\label{sldiffusion_model_cansys1}\\
    \pa_t q(x,t)&=2\ga P(x,t)+q(x,t)\left(\rho+b-\frac{2P(x,t)}{\left(1+P(x,t)^2\right)^2}\right)-D\Delta q(x,t),\label{sldiffusion_model_cansys2}\\
    \evalat{\pa_\nu P(x,t)}_{\pa\Omega}&=0,\quad
\evalat{\pa_\nu q(x,t)}_{\pa\Omega}=0,\label{sldiffusion_model_cansys3}\qquad 
    \evalat{P(x,t)}_{t=0}=P_0(x),\quad x\in\Omega. 
  \end{align}
\end{subequations}
We now explain how to use \poc\ to calculate CSS and canonical paths 
for \reff{slcan2}. 
For this we discuss files from the demo directory {\tt slocdemo} (except 
for obvious library files), assuming that {\tt slocdemo} is in the same 
directory as the libraries {\tt p2plib, p2poclib} and {\tt tom}. 
\subsubsection{Basics of \pdep, and the setup for CSS}\label{csssec}
We very briefly review the data structures of \pdep, and refer to 
\cite{p2p, p2p2} for more details and the underlying 
algorithms. The basic structure is a \mlab\ struct, henceforth called 
{\tt p} like {\tt p}roblem, which has a (large) number of fields 
(and subfields), as indicated in Table \ref{tab3}. 
\begin{table}[ht]
{\small
\begin{tabular}{p{17mm}|p{140mm}} 
field&purpose\\
\hline
fuha&struct of {\bf fu}nction {\bf ha}ndles; 
in particular the function handles 
p.fuha.sG, p.fuha.sGjac, p.fuha.bc, p.fuha.bcjac defining 
\reff{cs2a} and Jacobians.\\
nc, sw&{\bf n}umerical {\bf c}ontrols 
and {\bf sw}itches such as p.sw.bifcheck,\ldots \\[1mm]
u,np,nu& the solution u (including all parameters/auxiliary variables 
in u(p.nu+1:end)), the number of nodes p.np 
in the mesh, and the number of nodal values p.nu of PDE--variables\\
tau,branch&tangent tau(1:p.nu+p.nc.nq+1), and the branch, filled 
via bradat.m and p.fuha.outfu.\\
sol&other values/fields calculated at runtime, e.g.: ds (stepsize), res (residual), \ldots \\
usrlam&vector of user set target 
values for the primary parameter, default usrlam=[];\\ 
eqn,mesh&the tensors $c,a,b$ for the semilinear FEM setup, and 
the geometry data and mesh. \\
plot, file&switches (and, e.g., figure numbers and directory name) 
for plotting and file output\\
mat&problem matrices, e.g., mass/stiffness matrices $M$, $K$ for the 
the semilinear FEM setting. \\
\hline
\end{tabular}
\caption{{\small Selection (with focus on the semilinear case {\tt p.sw.sfem=1}) 
of fields in the structure {\tt p} describing a \pdep\ 
problem; see {\tt stanparam.m} in 
{\tt p2plib} for detailed information on the contents of these 
fields and the standard settings.} \label{tab3}}}
\end{table}
However, most of these can be set to standard values 
by calling {\tt p=stanparam(p)}. 
At least in simple problems, the user only has 
to provide: 
\bcen
\item The geometry of the domain $\Om$ and the boundary conditions.
\item Function handles (in the semilinear setting of interest here) 
{\tt sG} and, for speedup, {\tt sGjac}, implementing $G$, and its Jacobian. 
\item An initial guess for a solution $u$ of $G(u)=0$, i.e., an 
initial guess for a CSS. 
\ecen 

\begin{table}[!ht]
{\small
\begin{tabular}{|p{0.98\textwidth}|}\hline\\\vs{-10mm}
\begin{Verbatim}[commentchar=!, numbers=left] 
function p=slinit(p,lx,ly,nx,ny,sw) % init-routine 
p=stanparam(p); % set generic parameters to standard, if needed reset below..
p.nc.neq=2; p.fuha.sG=@slsG; p.fuha.sGjac=@slsGjac; %rhs 
p.fuha.outfu=@ocbra; p.fuha.jcf=@sljcf; p.fuha.con=@slcon; % current-val-obj 
p.usrlam=[0.55 0.6 0.65 0.7 0.75]; % target-values for bif-param lam 
[p.mesh.geo,bc]=recnbc2(lx,ly); p.vol=4*lx*ly; % geometry, and volume of dom 
p.fuha.bc=@(p,u) bc; p.fuha.bcjac=@(p,u) bc;   % standard Neumann BC 
p.xplot=lx; p.sw.spcalc=0; p.sw.jac=1; p.file.smod=100; % some more switches 
par=[0.03;0.55;0.5;0.5]; p.nc.ilam=2; % startup param values, and index of main param
% r=par(1); bp=par(2); cp=par(3); D=par(4); 
p.nc.dsmin=1e-6; p.nc.dsmax=0.5; p.nc.lammax=0.8; p.nc.lammin=0.549; p.sol.ds=0.1; 
p=stanmesh(p,nx,ny);p=setbmesh(p); p.sol.xi=1/p.np; % mesh 
p.eqn.c=[1;0;0;1;-1;0;0;-1]; p.eqn.a=0; p.eqn.b=0;  % diffusion tensor and a,b 
switch sw  % choose initial guess arcoording to switch 
case 1; u=0.3*ones(p.np,1); v=-13*ones(p.np,1); u0=[u v]; p.u=u0(:); % FSC
case 2; u=2*ones(p.np,1); v=-4*ones(p.np,1); u0=[u v]; p.u=u0(:); % FSM
case 3; .. % Scenario 2
end
p.u=[p.u; par]; p.sw.sfem=1; p=setfemops(p);  % semilin. setting 
[p.u,res]=nloop(p,p.u); fprintf('first res=%g\n',res); plotsol(p,1,1,1); 
\end{Verbatim}
\vs{-6mm}\\\hline\vs{-4mm}
\begin{Verbatim}[commentchar=!, numbers=left] 
function r=slsG(p,u) % CS for SLOC, p_t=D*lap p-1/q-b*p+p^2/(1+p^2)
%                                   q_t=-D lap q+2cp*p+q*(rho+bp-2*p/(1+p^2)^2; 
par=u(p.nu+1:end); r=par(1); bp=par(2); cp=par(3); D=par(4); 
P=u(1:p.np); q=u(p.np+1:2*p.np); 
f1=-1./q-bp*P+P.^2./(1+P.^2); f2=2*cp*P+q.*(r+bp-2*P./(1+P.^2).^2); 
f=[f1;f2]; r=D*p.mat.K*u(1:p.nu)-p.mat.M*f; 
\end{Verbatim}
\vs{-6mm}\\\hline\vs{-4mm}
\begin{Verbatim}[commentchar=!, numbers=left] 
function jc=sljcf(p,u) % current value J 
cp=u(p.nu+3:end); pv=u(1:p.np); kv=-1./u(p.np+1:p.nu); jc=log(kv)-cp*pv.^2;
\end{Verbatim} 
\vs{-6mm}\\\hline\vs{-4mm}
\begin{Verbatim}[commentchar=!, numbers=left] 
function k=slcon(p,u) % extract control from states/costates 
k=-1./u(p.np+1:p.nu); 
\end{Verbatim} 
\vs{-6mm}\\\hline
\end{tabular}}
\vs{0mm}
\caption{{\small The init routine {\tt slinit.m}, the rhs {\tt slsG.m}, 
the objective function {\tt sljcf.m}, and the function {\tt slcon.m}. See also, 
e.g.,  the source code 
of {\tt slsGjac} for the implementation of $G_u$. }\label{tab4}}\end{table}

Typically, the steps 1-3 are put into an init routine, 
here {\tt p=slinit(p,lx,ly,nx,ny,sw,rho)}, where {\tt lx,ly,nx,ny} 
are parameters to describe the domain size and discretization, 
and {\tt sw} is used to set up different initial guesses, 
see Table \ref{tab4}. The only 
additions/modifications to the standard \pdep\ setting for CSS problems 
are as follows: (the additional function handle) 
{\tt p.fuha.jc} should be set to the local current 
value objective function, here {\tt p.fuha.jc=@sljcf}, 
and {\tt p.fuha.outfu} to {\tt ocbra}, 
i.e., {\tt p.fuha.outfu=@ocbra}. This automatically puts $J_{ca}(u)$ 
at position 4 of the calculated output--branch. Here we generally use 
the {\em averaged} current objective function since typically 
we want to normalize $J_{ci}$ by the domain size for simple comparison 
between different domains. Finally, it is useful to set 
{\tt p.fuha.con=@slcon}, 
where {\tt k=slcon(p,u)} extracts the control $k$ from the states $v$, 
costates $\lam$ and parameters $\eta$, all contained in the vector 
{\tt u}.\footnote{Note that we do not use {\tt slcon} in {\tt slsG}. 
However, putting this function 
for the control into {\tt p} has the advantage that 
for instance plotting and extracting the value of the control can 
easily be done by calling some convenience functions of \poc.}

By calling {\tt p=cont(p)}, \pdep\ then first uses a Newton--loop to converge 
to a (numerical) solution, and afterwards 
attempts to continue in the given parameter. If {\tt p.sw.bifcheck>0}, 
then \pdep\ detects, localizes and saves to disk 
bifurcation points on the branch. Afterwards, the bifurcating branches 
can be computed by calling {\tt swibra} and {\tt cont} again. 
These (and other) \pdep\ commands (continuation, branch switching, and plotting) 
are typically put into a script file, here {\tt bdcmds.m}, 
see Table \ref{tab4b}.

\begin{figure}[ht]
\begin{tabular}{lll}
(a) BD of CSS&(b) BD, current values $J_{ca}$& (c) example CSS\\
\ig[width=0.3\textwidth]{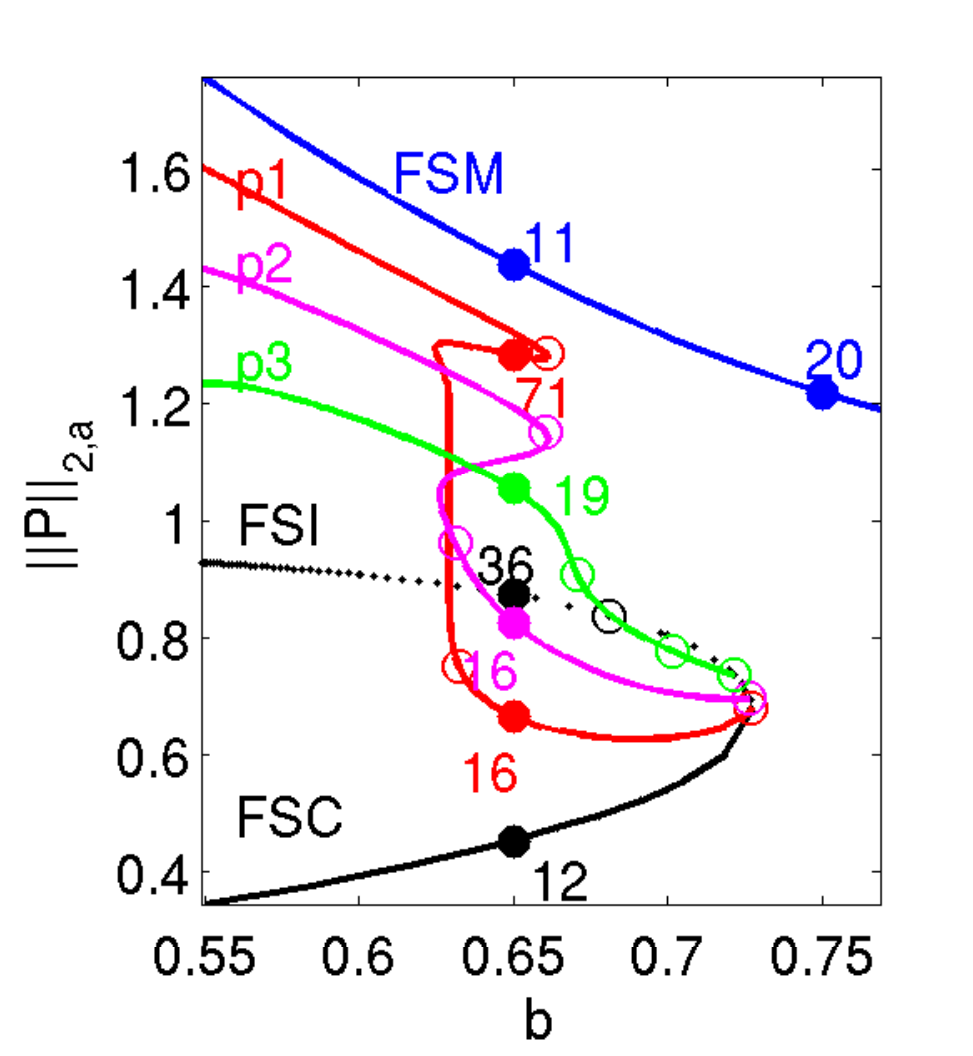}&
\ig[width=0.3\textwidth]{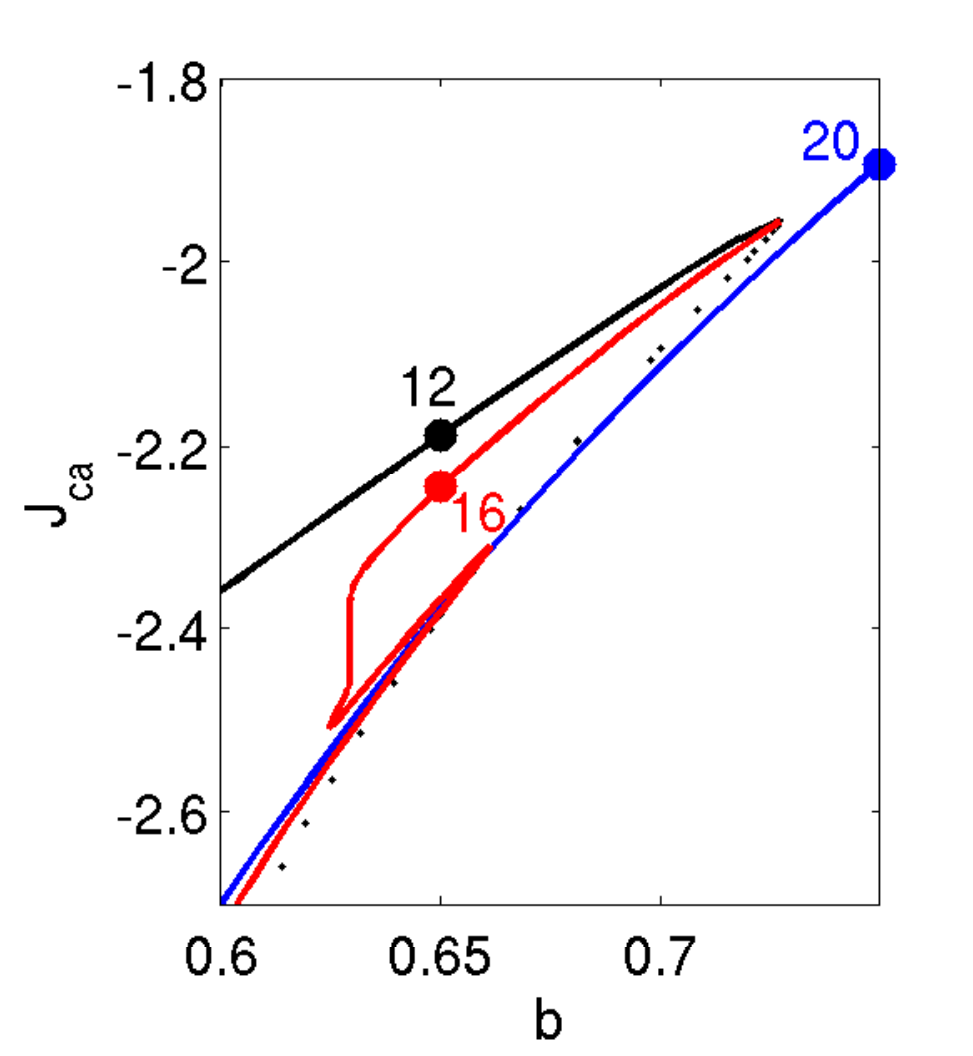}&
\raisebox{30mm}{\begin{tabular}{l}
\ig[width=0.15\textwidth]{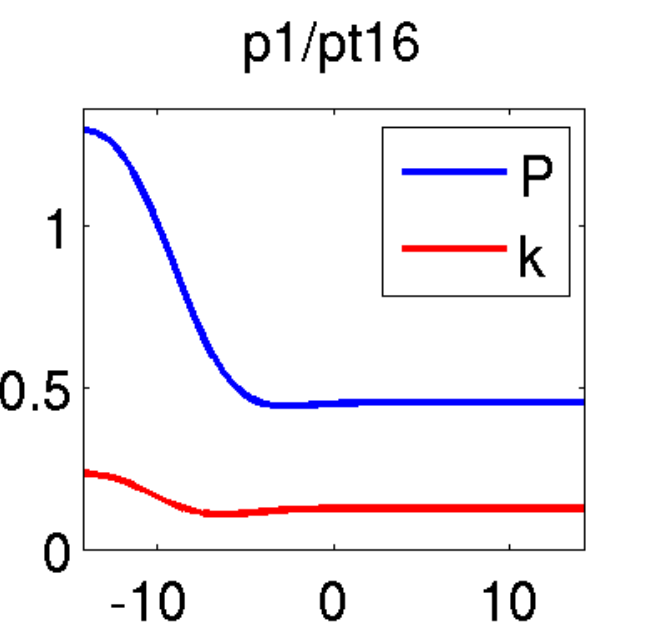}\ig[width=0.15\textwidth]{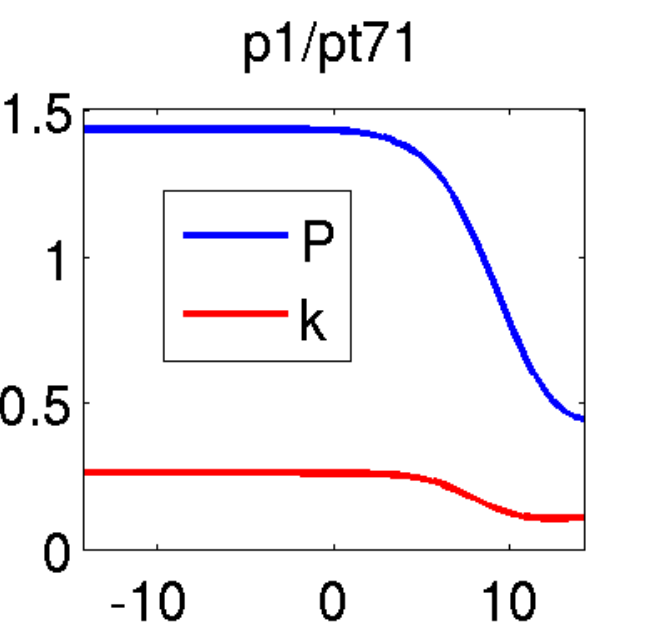}\\
\ig[width=0.15\textwidth]{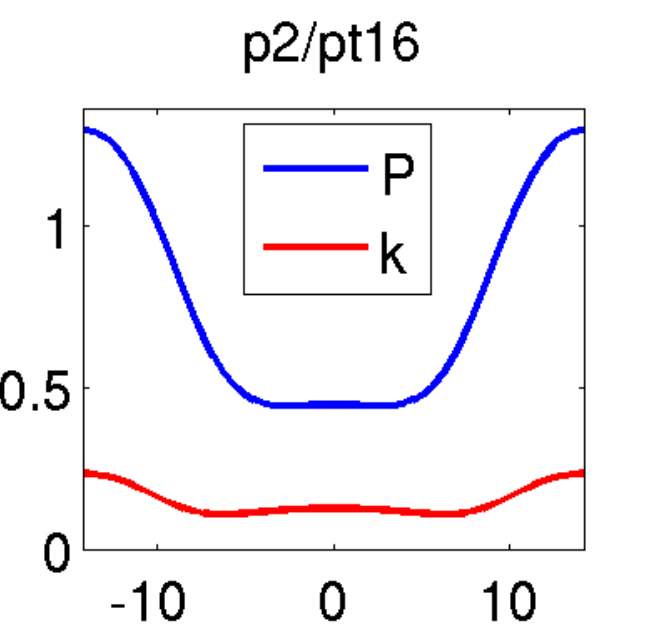}\ig[width=0.15\textwidth]{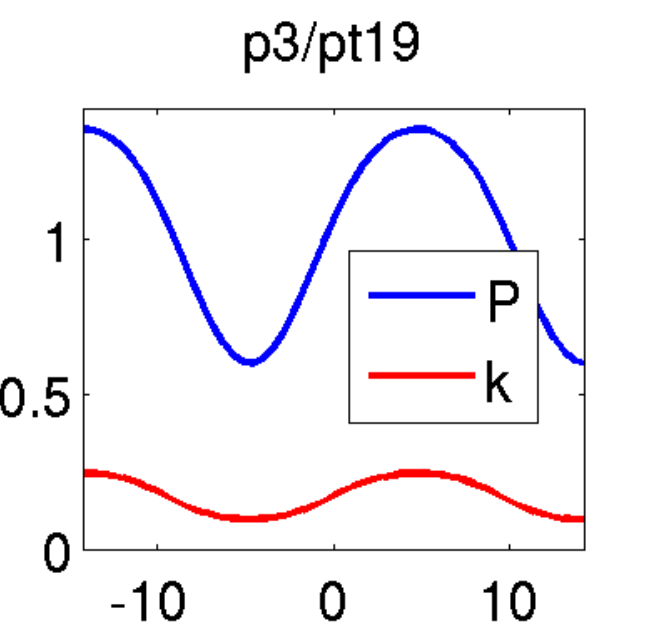}
\end{tabular}}
\end{tabular}
   \caption{{\small Example bifurcation diagrams and solution plots from running 
{\tt bdcmds.m}. For $b<b_{\text{fold}}\approx 0.73$ there are three branches of 
FCSS, here called FSC (Flat State Clean, low $P$), FSI (Flat State Intermediate), 
and FSM (Flat State Muddy, high $P$). On FSC$\cup$FSI there are a number 
of bifurcations to patterned CSS branches. The BD in $J_{ca}$ in (b) 
shows that at, e.g., $b=0.65$ we have $J(\text{FSI})<J(\text{FMS})<J(
\text{p1/pt71})<J(\text{p1/pt16})<J(\text{FSC})$, suggesting the FSC 
to be an optimal steady state. 
However, this does not exclude p1/pt16 or p1/pt71 to be 
optimal steady states, and moreover, there might be a path from the 
FSC to some other CSS with the SPP, dominating the FSC as a CSS. 
See \cite{GU15} for further discussion.} }
  \label{slfig1}
\end{figure}

\begin{table}[!ht]
{\small
\begin{tabular}{p{0.98\textwidth}}\hline\\[-8mm]
\begin{Verbatim}
%% stat. BD for sloc, main script-file. First set matlab paths: 
path('../p2plib',path);path('../p2poclib',path);
%% ---- Scenario 1, FSC/FSI branch 
close all; p=[];lx=2*pi/0.44; ly=0.1; nx=50; ny=1; sw=1; p=slinit(p,lx,ly,nx,ny,sw); 
p=setfn(p,'f1'); screenlayout(p); p=cont(p,100); 
%% ---- FSM branch 
sw=2; p=slinit(p,lx,ly,nx,ny,sw); p=setfn(p,'f2'); p.nc.dsmax=0.2; p=cont(p,15); 
%% ---- bif from f1 (set bpt* and p* and repeat as necessary) 
p=swibra('f1','bpt1','p1',-0.05); p.nc.dsmax=0.3; p.nc.neigdet=50; p=cont(p,150);
%% ---- plotting of BD, L2 and J_{ca}, and solution plots 
clf(3); pcmp=3; plotbraf('f1','bpt1',3,pcmp,'lab',12,'cl','k'); % FSC (+ other branches) 
clf(3); pcmp=4; plotbraf('f1','bpt1',4,pcmp,'lab',12,'cl','k'); % FSC (+other branches)
plot1Df('p1','pt16',1,1,1,2); plot1Df('p1','pt71',2,1,1,2); % solution plots ...
stancssvalf('p1','pt16'); % extract values <P>, <k>, J_{c,a} from solution in p1/pt16; 
\end{Verbatim}
\vs{-6mm}\\\hline
\end{tabular}}
\vs{0mm}
\caption{{\small Selected commands from the script file {\tt bdcmds.m}. 
See the source code for more details, and, e.g., {\tt bdcmds2D.m} 
for the quite similar commands in 2D. Note that these files are in \mlab\ 
cell-modes, and the cells should be executed one by one. 
In the last command, p.fuha.con must be set. 
}\label{tab4b}}\end{table}

Naturally,  
there are some modifications to the standard \pdep\ plotting commands, 
see, e.g., {\tt plot1D.m}. These work as usual by overloading the respective 
\pdep\ functions by putting the adapted file in the current directory. 
See Fig.~\ref{slfig1} 
for example results of running {\tt bdcmds}.

\subsubsection{Canonical paths}
The goal is to calculate canonical paths from some starting state $v(0)$ 
to a CSS $\uh_1$ with the SPP. For this we use one of the continuation 
algorithm {\tt iscnat} or {\tt iscarc} which in turn call {\tt mtom}, 
based on TOM. 
Since we only wanted minimal modifications of TOM we found it convenient 
(though somewhat dangerous) to pass a number of parameters to the 
functions called by \mtom\ via global variables. Thus, 
at the start of the canonical path scripts (here {\tt cpdemo.m}) we define 
a number of global variables, see Table \ref{tab5}. 

\begin{table}[ht]
{\small
\begin{tabular}{p{18mm}p{140mm}}
name&purpose\\
\hline
s0,s1&\pdep\ structs containing the boundary values 
at $t=0$ (s0) and at $t=T$ (s1)\\
Psi&the matrix $\Psi$ to encode the BC at $t=T$\\
u0, u1&vectors containing the current values of $u$ at the boundaries\\ 
par&the parameter values from s1 (only for convenience)\\
um1, um2&solutions at continuation steps $j-1$ and $j-2$ 
(to calculate secant predictors and used in extended system in {\tt iscarc})\\
sig&current (arclength) stepsize in {\tt iscarc}\\
\hline
\end{tabular}}
\caption{{\small Global variables for the computation of canonical paths, 
i.e., mainly for interfacing the 
driver scripts with the functions called by TOM.}\label{tab5}}
\end{table}
 
The usage of \poc\ to compute canonical paths 
is best understood by running and inspecting the 
demo file {\tt cpdemo.m}. 
\begin{figure}[ht]
\bce{\small
\begin{tabular}{ll}
(a) canonical path from p3/pt19 to FSC&
(b) A fold in $\al$\\
\ig[width=0.25\textwidth]{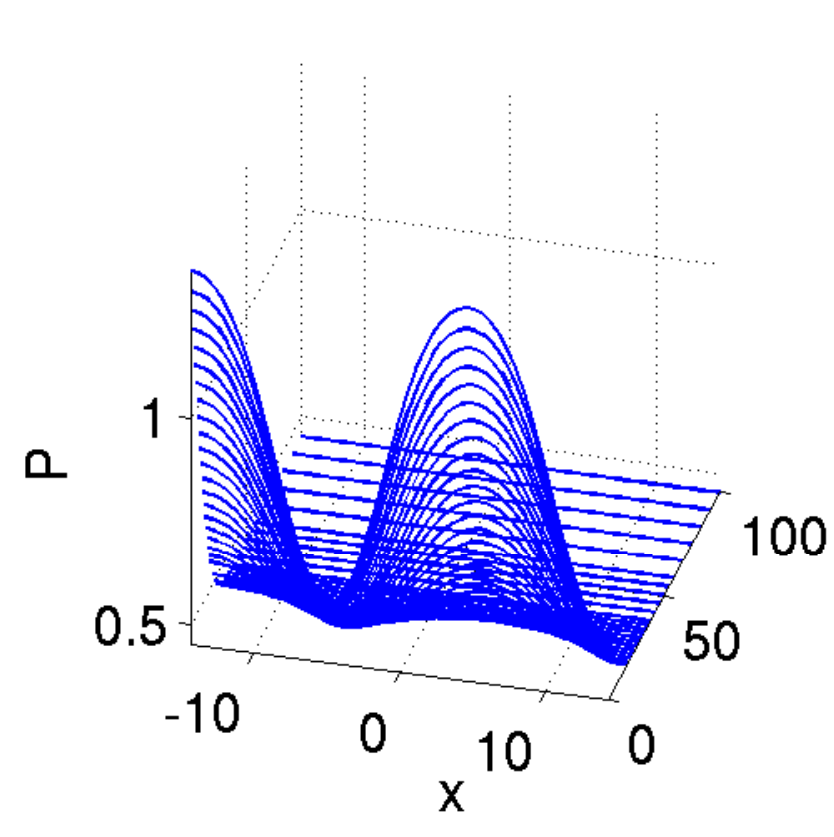}
\ig[width=0.25\textwidth]{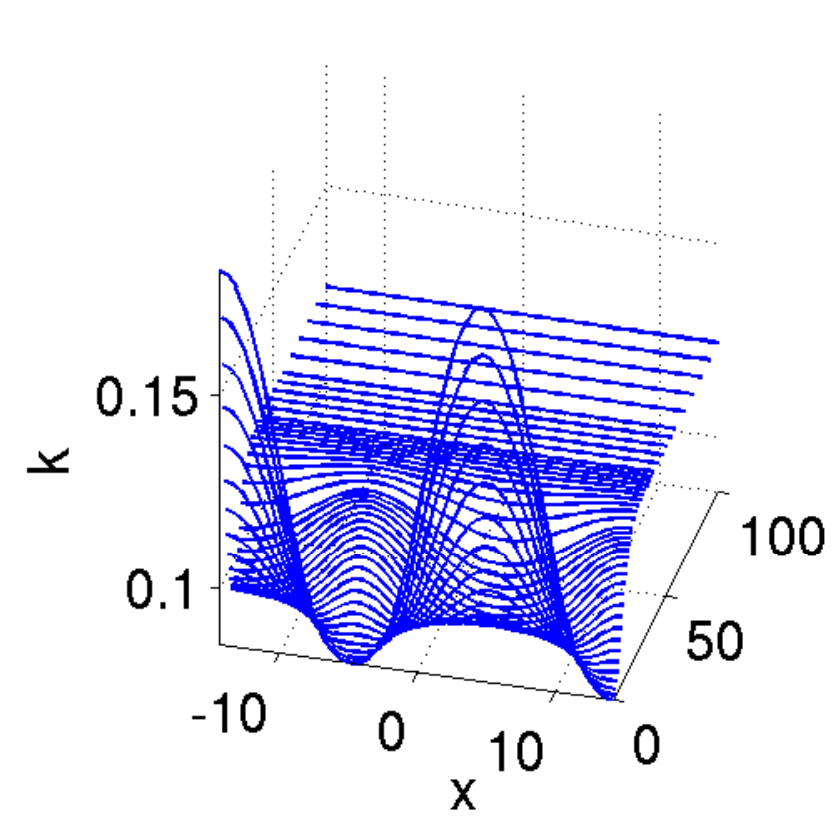}&
\ig[width=0.22\textwidth]{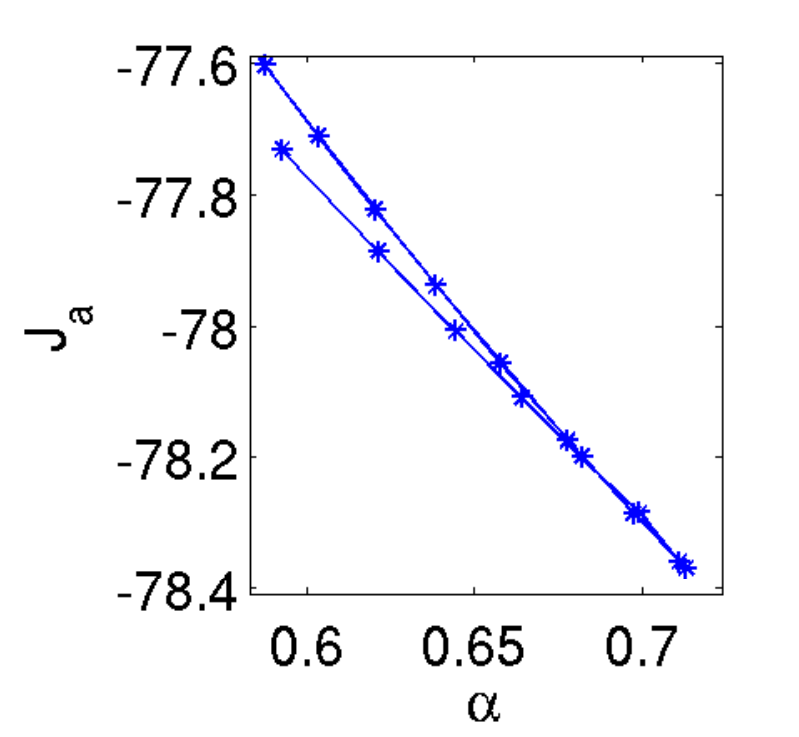}\\
(c) The ``upper'' canonical path at $\al=0.6$&(d) diagnostics\\
\ig[width=0.22\textwidth]{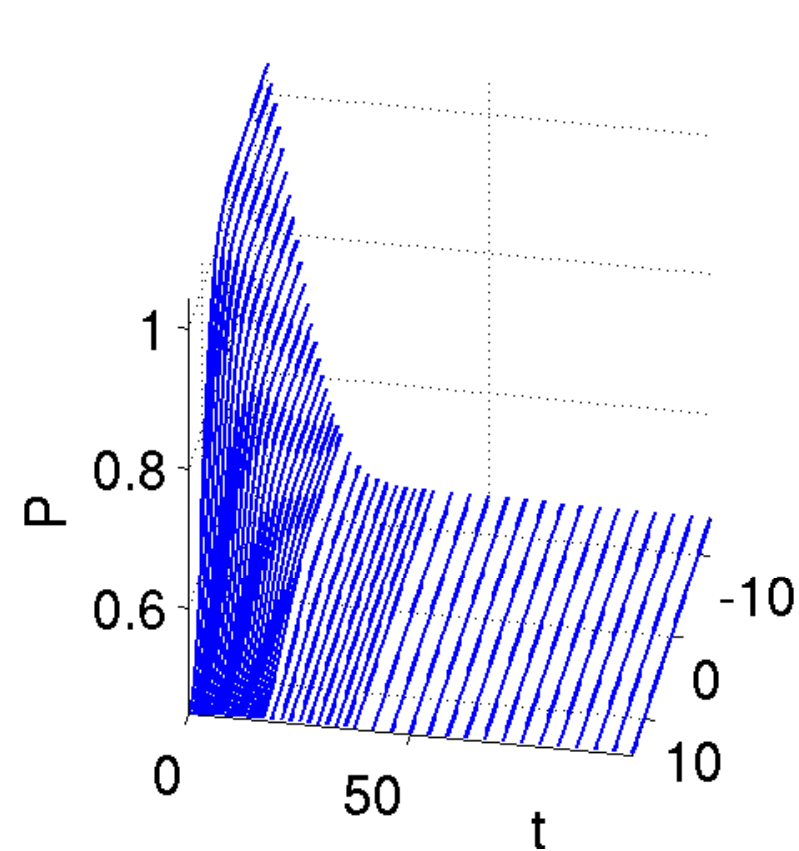}\ig[width=0.22\textwidth]{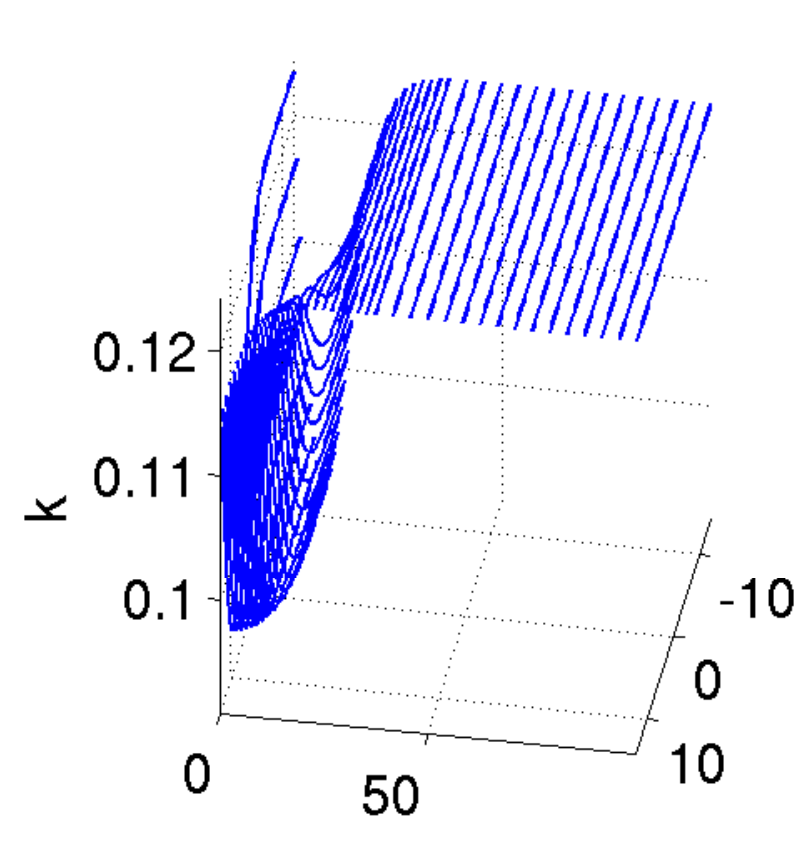}
&\raisebox{5mm}{\ig[width=0.2\textwidth]{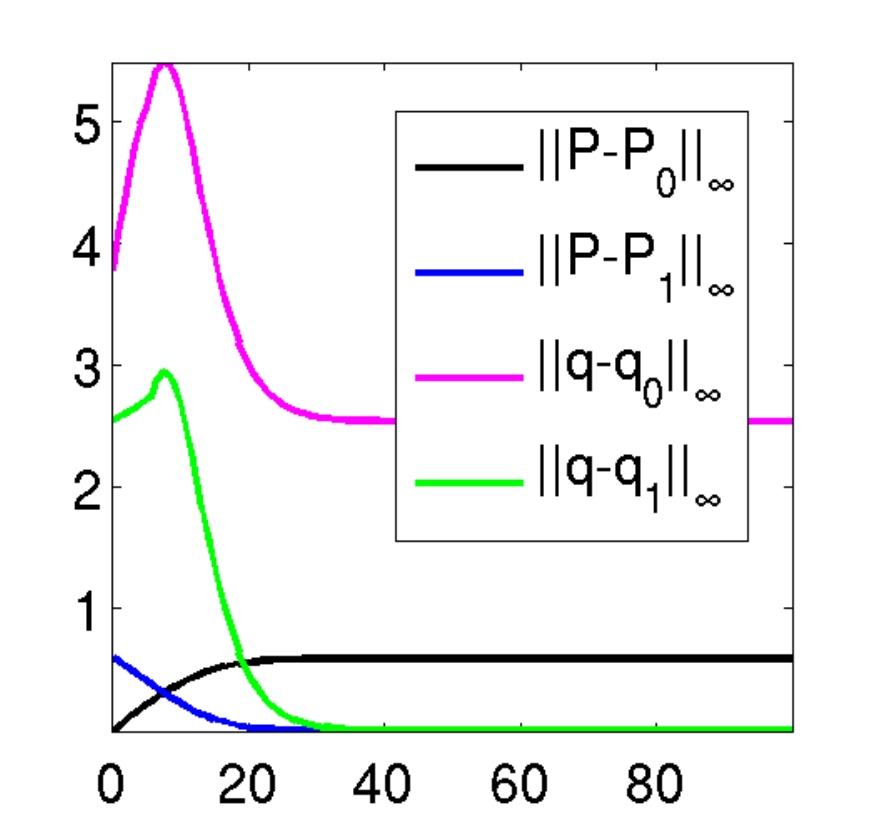}}\\
(e) The ``lower'' canonical path at $\al=0.6$&(f) diagnostics\\
\ig[width=0.22\textwidth]{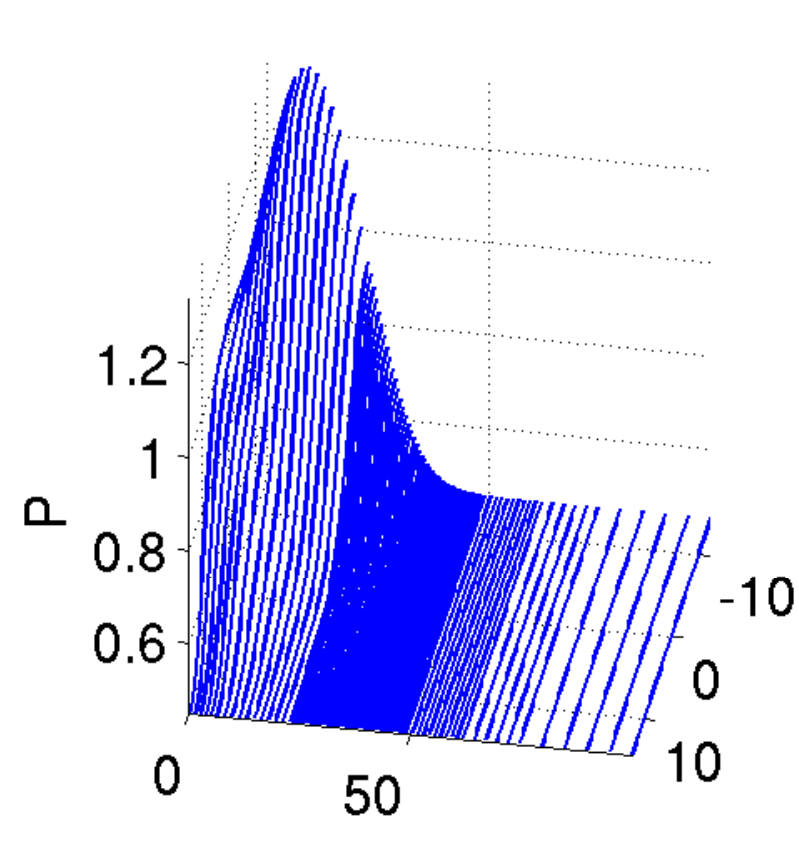}\ig[width=0.22\textwidth]{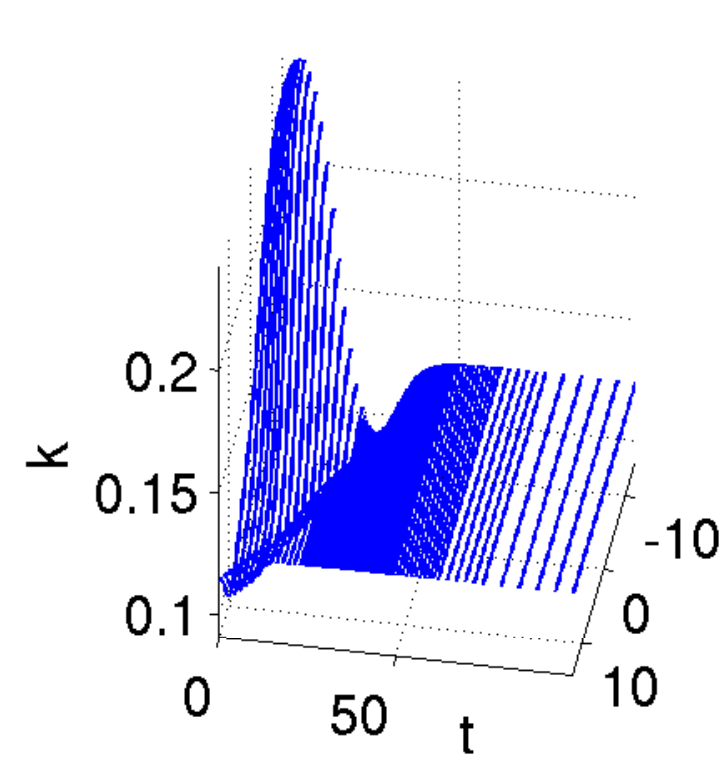}
&\raisebox{5mm}{\ig[width=0.2\textwidth]{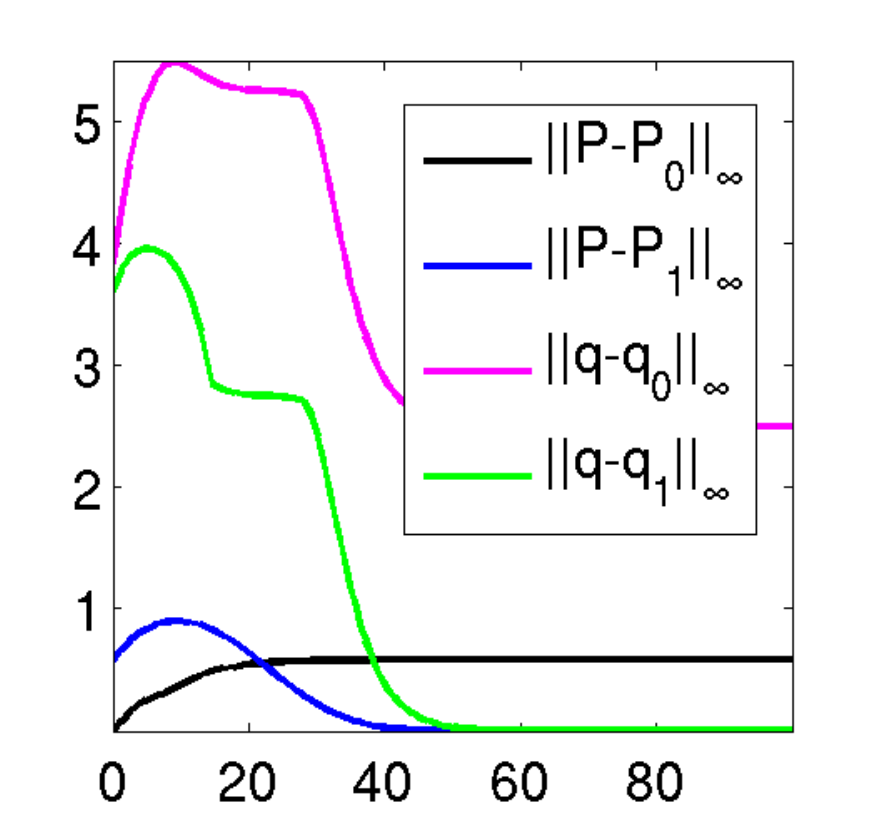}}
\end{tabular}
   \caption{Example outputs from running {\tt slcpdemo.m} \label{slfig2}}}
\ece
\end{figure}
Some results of running {\tt cpdemo} are shown in Fig.~\ref{slfig2}. (a) 
shows the ``easy'' case of a canonical path from  p3/pt19 to FSC 
(up to line 11), 
while (b) to (f) illustrate the case of a fold in $\al$ 
when trying to get a canonical path from p1/pt71 to the FSC. (c)-(e) show the 
two canonical paths obtained from picking two canonical paths from the output of {\tt iscarc} 
and a posteriori correcting to $\al=0.6$ (lines 24-26, only for the 
``upper'' canonical path). 
Line 21 from {\tt slcpdemo} 
prepares the calculation of Skiba paths, explained in \S\ref{skibasec}. 
We now give a brief overview of the involved \pocl\ functions, 
with the line-numbers refering to Table \ref{tab7}. 
\bci
\item[{\tt 
[alv,vv,sol,udat,tlv,tv,uv]=iscnat(alvin,sol,usec,opt,fn);}]  (line 9)\\ 
Input: 
{\tt alvin} as the vector of desired $\al$ values, for instance 
{\tt alvin=[0.25 0.5 1]}. {\tt sol,usec} 
can be empty (typically on first call), but on subsequent 
calls should contain the last solution  and 
the last secant (if {\tt opt.msw=1}). 
{\tt fn} is a structure containing the filenames for the start CSS and end CSS
\footnote{Taking $v$ from a CSS is basically for convenience: 
Of course, the initial states $v_0$ can be arbitrary, and 
there are no initial conditions for the co--states (and in particular 
those of the ``start CSS'' are not used). However, the construction is also 
motivated by the fact that one of the most interesting questions is if 
given $v$ from some CSS $\uh_0$ there exists a canonical path to 
an end CSS $\uh_1$, and whether this yields a higher value. 
Then, it is of course also interesting if one can also go the other 
way round, and for this we provide the {\tt flip} parameter in {\tt 
setfnflip}, see below.}, 
and {\tt opt} is an options structure containing TOM options 
and some more, see Table \ref{tab9} and Remark \ref{tomo}. \\
Output: {\tt alv} as the $\al$ vector of successful solves; 
{\tt vv} as the canonical path values $J_a$ of the successful solves; 
{\tt sol} as the  (last) 
canonical path; {\tt ydat} contains the 2 last steps and 
the last secant (useful for repeated calls, and for using 
{\tt iscnat} as startup for {\tt iscarc}). Finally, if {\tt opt.retsw=1}, 
then {\tt tlv,tv,uv} contain data of {\em all} successful solves, namely: 
for the $j$-th step, {\tt j=1:length(alv)}, 
{\tt tlv(j)}, contains the meshsize (in $t$),  
{\tt tv(j,1:tlv(j))} the mesh, and {\tt uv(j,1:n,1:tlv(j)} the solution. 
Thus, via {
{\tt sol.x=xv(j,1:tlv(j)); sol.y=squeeze(uv(j,\\1:n, 1:tlv(j)));} }  
the solution of the $j$-th step is recovered. This is useful 
for a posteriori inspecting some solution from the continuation 
(see line 24, and {\tt skiba.m} in {\tt vegdemo}).
Note that {\tt uv} can be large, and might give memory problems. 
If {\tt opt.retsw=0}, then {\tt tlv,tv,uv} are empty. 
\item[{\tt [alv,vv,usec,esol,tlv,tv,uv]=iscarc(esol,usec,opt,fn)};] (line 15)\\
Input: as for {\tt iscnat} (without {\tt alvin}), but with {\tt esol} 
containing the extended solution $(u_\al,\al)$, and similar for the secant 
{\tt usec}; as in {\tt iscnat}, if {\tt opt.start=1}, then {\tt esol,usec} 
can be empty. \\
Output: {\tt alv,val}$\in \R^{1\times m}$ as the vectors  
of achieved $\al$ and $J_a(\al)$; 
{\tt usec,esol} as the last secant/solution; {\tt tlv,tv,uv} as 
in {\tt iscnat},  but all in the sense of extended solutions, i.e. 
$(u_\al,\al)$. 
\eci 

\begin{table}
{\small
\begin{tabular}{p{15mm}p{60mm}|p{22mm}p{50mm}}
name&purpose&name&purpose\\
\hline
start&1 for startup, 0 else&M,lu,vsw&mass matrix, lu-switch and\\
retsw&1 to return full continuation data&&
 (extra) verbosity for mtom\\
rhoi&index of $\rho$ in par &t1&truncation time $T$\\
nti&\# of points in startup $t$--mesh&tv&current $t$--mesh\\
nsteps&number of steps for {\tt iscarc}&sigmin,sigmax&
min\&max stepsize for {\tt iscarc}\\
\hline
\end{tabular}}
\caption{{\small Switches/controls in {\tt opt} besides the TOM options. 
Note that {\tt iscarc} only 
has a very elementary stepsize control via {\tt opt.sigmin,opt.sigmax}. 
\label{tab9}}}
\end{table}

\brem{\label{tomo}\rm 
Concerning the original TOM options we remark that typically we run 
{\tt iscnat} and {\tt iscarc} with weak error requirements and what appears to be the 
fastest monitor and order options, i.e., 
{\tt opt.Monitor=3; opt.order=2;}. Once continuation is successful 
(or also if it fails at some $\al$), we can always postprocess by calling 
{\tt mtom}  
again with a higher order, stronger error requirements, and different 
monitor options, e.g., mesh--refinement based on condition rather than error 
alone. See the original TOM documentation. 
}\erem 

\begin{table}[!ht]
{\small
\begin{tabular}{p{0.98\textwidth}}\hline\\\vs{-8mm}
\begin{Verbatim}[commentchar=!, numbers=left] 
% driver script for Shallow Lake Optimal Control, first set paths and globals 
path('../tom',path); path('../p2poclib',path); 
close all; clear all; global s0 s1 u0 u1 Psi par xi ym1 ym2 sig;  
%% Preparations: put filenames into fn, set some bvp parameters 
sd0='f1'; sp0='pt12'; sd1='p3'; sp1='pt19'; flip=1; % p3->FSC
fn=setfnflip(sd0,sp0,sd1,sp1,flip); opt=[]; opt=ocstanopt(opt);
opt.rhoi=1; opt.t1=100; opt.start=1; opt.tv=[]; opt.nti=10; opt.retsw=0; 
%% the solve and continue call, and some plots 
sol=[]; alvin=[0.1 0.25 0.5 0.75 1]; v=[15,30];
[alv,vv,sol,ydat,tlv,xv,yv]=iscnat(alvin,sol,[],opt,fn); slsolplot(sol,v); 
%% ---- A fold in alpha, here iscarc needed. Prep. and initial iscarc call
sd0='f1'; sp0='pt12'; sd1='p1'; sp1='pt71'; flip=1; fn=setfnflip(sd0,sp0,sd1,sp1,flip); 
esol=[]; ysec=[]; opt.nsteps=3; opt.alvin=[0.2 0.25]; sig=0.1; opt.nti=10; opt.tv=[];
opt.Stats_step='on'; opt.start=1; opt.sigmax=1; opt.retsw=1; 
[alv,vv,ysec1,esol1,tlv,xv,yv]=iscarc(esol,ysec,opt,fn); opt.start=0; 
%% subsequent iccarc-calls (repeat this cell) 
opt.nsteps=20; ysec=ysec1; esol=esol1; % new input (for repeated calls) 
[alv1,vv1,ysec1,esol1,tlv1,xv1,yv1]=iscarc(esol,ysec,opt,fn); 
alv=[alv alv1]; vv=[vv vv1]; tlv=[tlv tlv1]; xv=[xv; xv1]; yv=[yv; yv1]; 
%% Postprocess sol from iscarc, first a simple plot of J over alpha  
alv0=alv; vv0=vv; xv0=xv; yv0=yv; tlv0=tlv; % save results for skibademo.m 
figure(6); clf; plot(alv(1,:),vv(1,:),'-*');  xlabel('\alpha'); ylabel('J_{a}');
%% fix al from iscarc to some given value and compute CPs, first j=22, then j=34
j=22; tl=tlv(j); n=s1.nu; sol.x=xv(j,1:tlv(j));sol.y=squeeze(yv(j,1:n,1:tlv(j))); 
al=0.6; u0=al*s0.u(1:n)+(1-al)*s1.u(1:n); u1=s1.u(1:s1.nu); 
opt.M=s1.mat.M; sol=mtom(@mrhs,@cbcf,sol,opt); v=[100,30]; slsolplot(sol,v); 
\end{Verbatim}
\vs{-6mm}\\\hline
\end{tabular}}
\vs{0mm}
\caption{{\small Selected commands from the script file {\tt cpdemo.m}. 
See the source code for more details and, e.g., more fancy 
plotting. }\label{tab7}}\end{table} 
The functions {\tt iscnat} and {\tt iscarc} are the two main user 
interface functions for the canonical path numerics. There 
are a number of additional functions for internal use, and some convenience 
functions, which we briefly review as follows:

\bci 
\item[{\tt [Psi,mu,d,t]=getPsi(s1)};] compute $\Psi$, the eigenvalues {\tt mu}, 
the defect d, and a suggestion for $T$. Note that this becomes expensive 
with large $2nN$ (i.e., the total number of DoF). 
\item[{\tt [sol,info]=mtom(ODE,BC,solinit,opt,varargin)};] 
the ad--hoc modification of TOM, which allows for $M$ in (\ref{bvp1}a). 
Extra arguments $M$ and {\tt lu,vsw} in {\tt opt}. 
If opt.lu=0, then $\backslash$ is used for solving linear systems 
instead of an LU--decomposition, 
which becomes too slow when $nN\times m$ 
becomes too large. See the TOM documentation for all other arguments 
including {\tt opt}, 
and note that the modifications in mtom can be identified by searching 
``HU'' in  {\tt mtom.m}. Of course {\tt mtom} (as any other 
function) can also be called directly (line 26), which for instance is useful 
to postprocess the output of some continuation by changing parameters by hand. 
\item[{\tt f=mrhs(t,u,k); J=fjac(t,u);} and {\tt  f=mrhse(t,u,k); J=fjace(t,u);}]
the rhs and its Jacobian to be called within {\tt mtom}. These are just wrappers 
which calculate $f$ and $J$ by calling the resp.~functions in the 
\pdep--struct {\tt s1}, which were already set up and used to calculate 
the CSS. {\tt s1} is passed as a global variable. 
Similar remarks apply to {\tt mrhse} and {\tt fjace} 
for the extended setting in {\tt iscarc}. 
\item[{\tt bc=cbcf(ya,yb);[ja,jb]=cbcjac(ya,yb)} and 
{\tt bc=cbcfe(ya,yb);[ja,jb]=cbcjace(ya,yb)};] The boundary conditions 
(in time) for \reff{bvp1} and the associated Jacobians. Implemented by passing 
$u_0, \uh_1, \Psi$ and similar globally. The {\tt *e} (as in {\tt e}xtended) 
versions are 
for {\tt iscarc}. 
\item[{\tt jcaval=jcai(s1,sol,rho)} and {\tt djca=isjca(s1,sol,rho);}] 
Calculate the objective value 
\huga{\label{jadef} 
J(u)=\int_0^T \er^{-\rho t}J_{ca}(v(t,\cdot),k(t,\cdot))\dd t}
of the solution $u$ in {\tt sol} (with $J_c$ taken from {\tt s1.fuha.jcf}), 
and similarly the normalized discounted value of a CSS contained 
in {\tt sol.y(:,end)}. 
\item[{\tt fn=setfnflip(sd0,sp0,sd1,sp1,flip);}] generate the filename struct 
{\tt fn} from sd0, sp0 (sd0/sp0.mat contains IC $u_0$) and sd1,sp1 
(contains $\uh$); if flip=1, then interchange *0 and *1. 
\item[{\tt psol3D(p,sol,wnr,cmp,v,tit);}] $x$--$t$ plots of canonical paths;  
plot component {\tt cmp} of a canonical path 
{\tt sol} to figure {\tt wnr}, with view v and title {\tt tit}. 
If {\tt cmp=0}, then plot the control {\tt k}, extracted from {\tt sol} 
via {\tt p.fuha.con}.
\eci 
Thus, after having set up {\tt p} as in \S\ref{csssec} for the CSS, 
including $G$ and {\tt p.fuha.jcf}, the user does not need to set up 
any additional functions to calculate canonical paths and their values. 
However, typically there are some functions which 
should be adapted to the given problem, e.g., for plotting, for instance 
\bci 
\item[{\tt slsolplot(sol,v);}] (line 10) which calls: 
\item[{\tt zdia=sldiagn(sol,wnr);}] Plot some norms on a canonical path 
as functions of $t$ to figure(wnr). This is for instance useful 
to check the convergence behaviour of the canonical path as $t\ra T$, 
cf.~(d),(f) in Fig.~\ref{slfig2}. 
\eci

\subsubsection{A patterned Skiba point}\label{skibasec}
In ODE OC applications, if there are several locally stable OSS, 
then often an important issue is to identify their domains of attractions. 
These are separated by so called theshhold or Skiba--points (if $N=1$) 
or Skiba--manifolds (if $N>1$), see \cite{Skiba78} and \cite[Chapter 5]{grassetal2008}. 
Roughly speaking, these are initial states from which there are several 
optimal paths with the same value but leading to different CSS. In PDE 
applications, even under spatial discretization with moderate $nN$, 
Skiba manifolds should be expected to become very complicated objects. 
Thus, here we just give one example how to compute a patterned Skiba 
point between FSC and FSM. 

\begin{figure}[ht]
\bce{\small
\begin{tabular}{ll}

(a) A Skiba point at $\al=0.454$&(b) Paths to FSC (blue) and FSM (red)\\
\ig[width=0.28\textwidth,height=35mm]{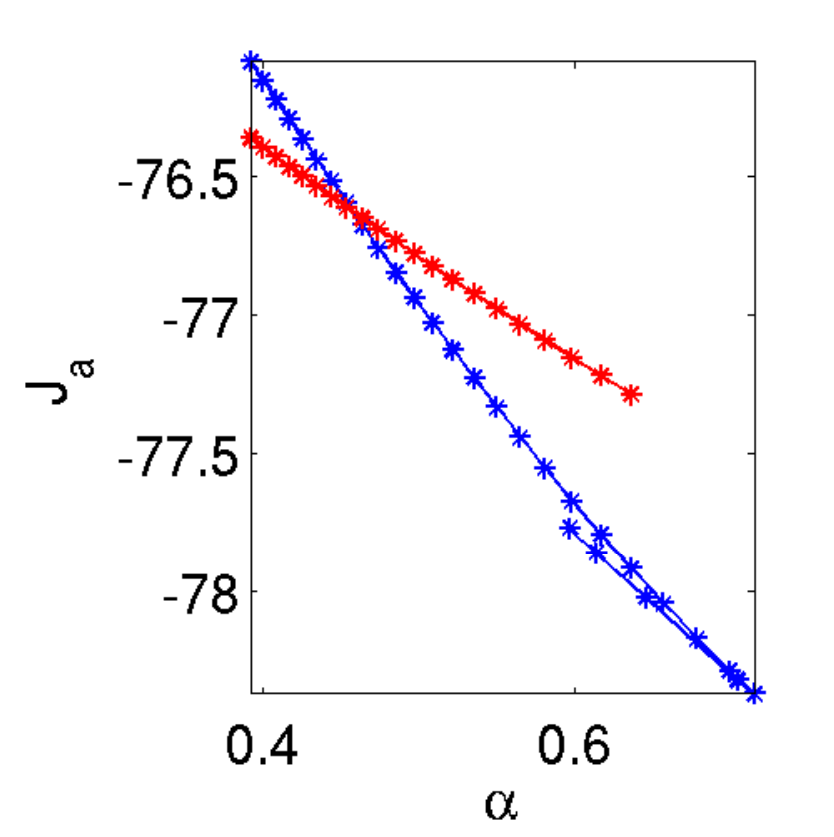}&\ig[width=0.22\textwidth]{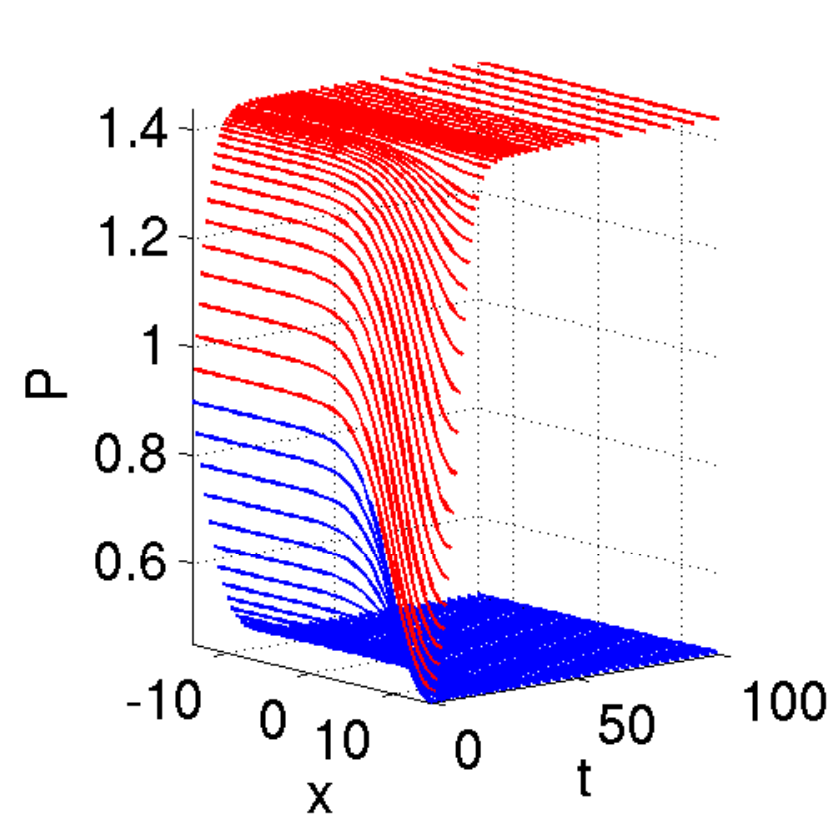}
\raisebox{0mm}{\ig[width=0.2\textwidth]{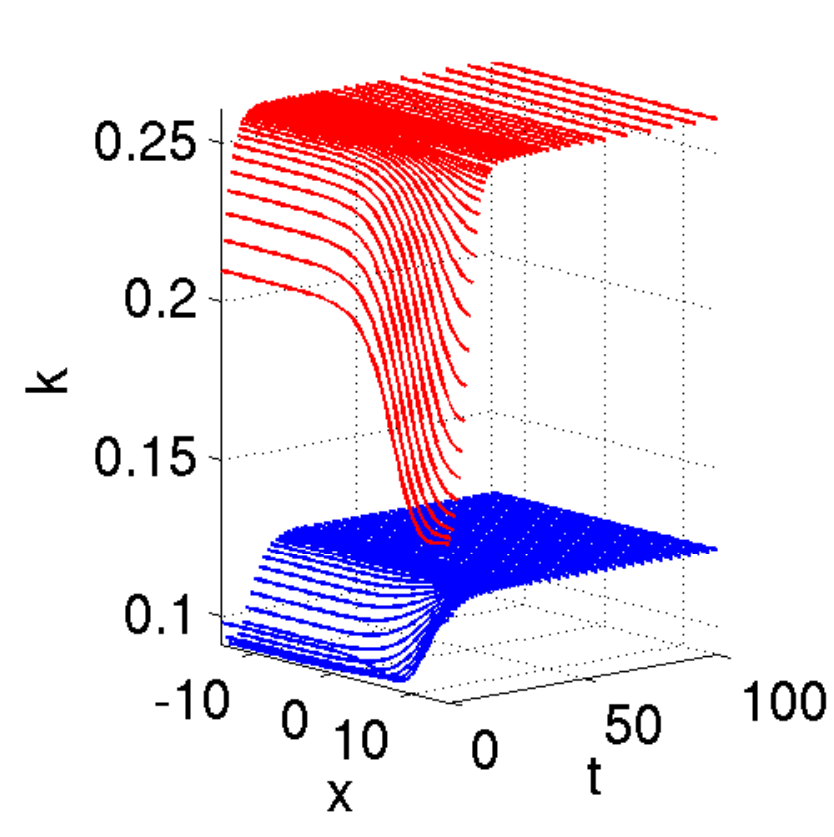}}
\end{tabular}
   \caption{Example outputs from {\tt skibademo.m} \label{slfig3}}}
\ece
\end{figure}

In Line 17-19 of {\tt cpdemo.m} we attempt to find a path from $P_{\text{PS}}$ 
given by 
{\tt p1/pt71} to $(P,q)_{\text{FSC}}$ given by {\tt FSC/pt12}; this fails 
due to the fold in 
$\al$. However, for given $\al$ we can also try to find a path from 
the initial state 
$P_\al(0):=\al P_{\text{PS}}+(1-\al)P_{\text{FSC}}$ to the FSM, and compare to 
the path to the FSC. For this, in line 21 of {\tt cpdemo.m}, we stored
the $\al$ and $J_{ca}$ values into {\tt alv0, vv0}, 
and also the path data into {\tt tlv0,tv0 uv0}. See {\tt skibademo.m} 
in Table \ref{tab7b} 
(in particular line 11 and the following) how to put the 
values {\tt uv0(j,:,1)} into {\tt s0} and 
subsequently find the paths to the FSM, and Fig.~\ref{slfig3} for 
illustration. 

\begin{table}[!ht]
{\small
\begin{tabular}{p{0.98\textwidth}}\hline\\\vs{-8mm}
\begin{Verbatim}[commentchar=!, numbers=left]
% Skiba example, continues cpdemo.m 
% find paths from the yv0 initial states from cpdemo.m to FSM 
js=10; je=30; jl=js-je+1; % alpha-range; now set the target and Psi to FSM: 
s1=loadp('f2','pt11'); u1=s1.u(1:s1.nu); [Psi,muv,d,t1]=getPsi(s1); 
a0l=length(alv0); tva=zeros(jl,opt.Nmax+1); % some prep. and fields to hold paths 
uva=zeros(jl,n+1,opt.Nmax+1); alva=[]; vva=[]; tavl=[]; sol=[]; 
alvin=[0.1 0.25 0.5 0.75 1]; % we run from uv0(j,:) to FSM with iscnat
tv=linspace(0,opt.t1,opt.nti); se=2; opt.tv=tv.^se./opt.t1^(se-1); doplot=1; 
opt.msw=0; opt.Stats_step='off'; v=[50,8]; % switch off stats 
for j=js:je; 
   fprintf('j=%i, al=%g\n', j, alv0(j)); s0.u(1:n)=uv0(j,1:n,1)'; 
   [alv,vv,sol,udat]=iscnat(alvin,[],[],opt,fn); 
   if alv(end)==1; Jd=vv0(j)-vv(end); fprintf('J1-J2=%g\n',Jd); % contin. successful 
      alva=[alva alv0(j)]; vva=[vva vv(end)]; tl=length(sol.x); % put vals in vector
      tavl=[tavl tl]; tva(j,1:tl)=sol.x; uva(j,1:n,1:tl)=sol.y; 
      if abs(Jd)<0.05; doplot=asknu('plot path?',doplot); % Skiba point(s) found
        if doplot==1; sol0=[]; alp=alv0(j); % plot the paths to FSC and FSM 
          sol0.x=tv0(j,1:tlv0(j));sol0.y=squeeze(uv0(j,1:n,1:tlv0(j))); 
          psol3Dm(s1,sol0,sol,1,1,[]); view(v); zlabel('P'); pause
        end
      end
   end
end
%% plot value diagram 
figure(6); plot(alv0(js:jep),vv0(js:jep),'-*b');hold on;plot(alva,vva,'-*r');
xlabel('\alpha','FontSize',s1.plot.fs); ylabel('J_{a}','FontSize',s1.plot.fs);
\end{Verbatim}
\vs{-6mm}\\\hline
\end{tabular}}
\vs{0mm}
\caption{{\small The script file {\tt skibademo.m}. See text for comments. 
}\label{tab7b}}\end{table} 

\subsubsection{Further comments}
To keep the demos simple, here we do {\em not} include 
versions of {\tt iscnat} and {\tt iscarc} that use the truncation time 
$T$ as an additional free parameter \cite[Fig.5]{GU15}. However, 
the files {\tt bdcmds.m} and {\tt cpdemo.m} also contain some of 
the commands used to study Scenario 2, see \cite[\S2.3]{GU15} for the results, 
and the directory {\tt slocdemo} contains the script files 
{\tt bdcmds2D.m} and {\tt cpdemo2D.m}, used to compute CSS and canonical 
paths for \reff{slcan2} over the 2D domain 
$\Om=(-L,L)\times (-\frac L 2, \frac L 2)$ (based on exactly the same 
init file {\tt slinit.m}), some auxiliary plotting functions 
{\tt plotsolf.m} and {\tt plotsolfu.m}, and the function {\tt sol2mov.m} used 
to generate movies of canonical paths. See, e.g., \cite[Fig.~7,8]{GU15} 
for some results, while 
some movies can be downloaded at the \pdep\ homepage.


\subsection{The vegOC model}
Our second example, from \cite{U15}, 
concerns the optimal control of a reaction diffusion 
{\em system} used to model grazing in a semi arid system for 
biomass (vegetation) $v$ and soil water $w$, following \cite{BX10}. 
Here, semi arid means that there is enough water to support {\em some} 
vegetation, but not enough water for a dense homogeneous vegetation. 
This is an important problem as it is estimated that semi arid areas 
cover about 40\% of the world's land area and support about 
two billion people, often by grazing livestock, \url{www.allcountries.org/maps/world_climate_maps.html}. In semi arid areas, often overgrazing is a serious 
threat as it may lead to irreversible desertification, see, e.g., \cite{SBB09}, 
and the references therein. 

Denoting the harvesting (grazing) effort as the control by $E$, we consider 
\begin{subequations}\label{v1}
\hual{V(v_0,w_0)&=\max_{E(\cdot,\cdot)} J(v_0,w_0,E), \\
\pa_t v&=d_1\Delta v+[gwp^\eta-d(1+\del v)]v-H,\\
\pa_t w&=d_2\Delta w+R(\beta+\xi v)-(r_u v+r_w)w, 
}
with harvest $H=v^\al E^{1-\al}$, and current value 
objective function $J_{c}=J_c(v,E)=pH-cE$, which thus depends on the price $p$, 
the costs $c$ for harvesting/grazing, and $v$, $E$ in a classical 
Cobb--Douglas form 
with elasticity parameter $0<\al<1$. For the modeling, and the meaning and 
values of the parameters 
$(g,\eta,d,\delta,\beta,\xi, r_u,r_w, d_{1,2})$ we refer to 
\cite{BX10,U15} and the references therein (see also 
Table \ref{vtab1}, line 9 for the parameter values), and here only remark that 
$\rho=0.03, p=1.1, \al=0.3, c=1$ are the economic parameters, and we take 
the rainfall $R$ as the main 
bifurcation parameter. Furthermore, we have the BC and IC 
\huga{\label{nbc}
\pa_\nu v=\pa_\nu w=0 \text{ on } \pa \Om,\quad (v,w)|_{t=0}=(v_0,w_0).  
}
\end{subequations}
Denoting the co-states by $(\lam,\mu)$ we have the local current value Hamiltonian 
\hual{
 \CH(v,w,\lam,\mu,E)=J_c(v,E) 
&+\lam\bigl[d_1\Delta v+(gwp^\eta-d(1+\del v))v-H\bigr]\notag\\
&+\mu\bigl[d_2\Delta w+R(\beta+\xi v)-(r_u v+r_w)w\bigr], 
} 
and obtain the canonical system 
\begin{subequations}\label{vegcs}
\hual{
\pa_t v&=\CH_\lam=d_1\Delta v+[gwp^\eta-d(1+\del v)]v-H,\\
\pa_t w&=\CH_\mu=d_2\Delta w+R(\beta+\xi v)-(r_u v+r_w)w,\\
\pa_t\lam&=\rho\lam-\CH_v=\rho\lam-p\al v^{\al-1}E^{1-\al}-
\lam\bigl[g(\eta+1)wv^\eta-2d\del v-d-\al v^{\al-1}E^{1-\al}]\\\notag 
&\qquad\qquad\qquad-\mu(R\xi-r_u)w-d_1\Delta \lam, 
\\
\pa_t\mu&=\rho\mu-\CH_w=\rho\mu-\lam g v^{\eta+1}+\mu(r_u v+r_w)-d_2\Delta \mu, 
}
where $E$ is obtained from solving $\pa_E \CH=0$, giving 
\huga{
E=\left(\frac{c}{(p-\lam)(1-\al)}\right)^{-1/\al}v. 
}
With the notation $u=(v,w,\lam,\mu)$, the IC, the BC, and the 
transversality condition are 
\huga{\label{nbc2}
(v,w)|_{t=0}=(v_0,w_0), \quad 
\pa_\nu u=0 
\text{ on } \pa \Om, \quad \lim_{t\to\infty}\er^{-\rho t}u(t)=0. 
}
\end{subequations}

To study \reff{vegcs}, we write it as $\pa_t u=-G(u)$ and basically need to
 set up $G$ and the BC. This follows the general \pdep\ settings with 
the OC related modifications already explained in \S\ref{slsec}, and 
thus we only give the following remarks, first concerning {\tt veginit.m}, 
see Table \ref{vtab1}. 
\bci 
\item[In line 2] we only set up {\tt p.fuha.sG} since in this demo we 
use {\tt p.sw.jac=0} (numerical Jacobians), 
and hence do not need to set {\tt p.fuha.sGjac}. 
\item[lines 4-7] set the Neumann BC and diffusion tensor for the 4 component 
system (see {\tt gnbc.m} and {\tt isoc.m} for documentation) 
\item[lines 8-10] set the desired $R$ values for output of CSS to disk, 
the parameter values, and the main bifurcation parameter. 
Of course, one could also hard-code all parameters except $R$, but 
we generally recommend to treat parameters as parameters 
since this is needed if later a continuation in some other parameter 
is desired, and since it usually makes the code more readable. 
\eci 
\begin{table}[!ht]
{\small
\begin{tabular}{p{0.98\textwidth}}\hline\\\vs{-8mm}
\begin{Verbatim}[commentchar=!, numbers=left] 
function p=veginit(p,lx,ly,nx,ny,sw,rho) % init-routine for vegOC 
p=stanparam(p); p.nc.neq=4; p.fuha.sG=@vegsG; p.fuha.jcf=@vegjcf; p.fuha.outfu=@ocbra;
p.mesh.geo=rec(lx,ly); p=stanmesh(p,nx,ny); p.sol.xi=0.005/p.np; % generate mesh
q=zeros(p.nc.neq); g=zeros(p.nc.neq,1); % setting up Neumann BC for 4 components 
bc=gnbc(p.nc.neq,4,q,g); p.fuha.bc=@(p,u) bc; p.fuha.bcjac=@(p,u) bc;
p.d1=0.05; p.d2=10; p.eqn.a=0; p.eqn.b=0; % setting up K for 4 components 
c=diag([p.d1, p.d2, -p.d1, -p.d2]); p.eqn.c=isoc(c,4,1); 
p.usrlam=[4 10 20 26 28]; % desired R values for output of CSS 
par=[rho 1e-3 0.5 0.03 0.005 0.9 1e-3 34 0.01 0.1 1 1.1 0.3]; % par-values 
p.nc.ilam=8;  % choose the active par, here Rainfall R 
% now continue with setting a few more param and the initial guess ... see veginit.m
\end{Verbatim}
\vs{-6mm}\\\hline
\end{tabular}}
\vs{0mm}
\caption{{\small Selected commands from the init-routine {\tt veginit.m}. 
See the source code for more details. }\label{vtab1}}\end{table}
Table \ref{vtab2} shows the complete codes for setting up $G$ and the 
current value $J_c$. Both use the auxiliary function {\tt efu}, 
which is also used in, e.g., {\tt valf.m} to tabulate characteristical values 
of CSS.  

\begin{table}[!ht]
{\small
\begin{tabular}{p{0.98\textwidth}}\hline\\\vs{-8mm}
\begin{Verbatim}[commentchar=!, numbers=left] 
function r=vegsG(p,u) % rhs for vegOC problem 
par=u(p.nu+1:end); rho=par(1); g=par(2); eta=par(3); % extract param 
d=par(4); del=par(5); beta=par(6); xi=par(7); rp=par(8); up=par(9); rw=par(10); 
cp=par(11); pp=par(12); al=par(13); [e,h,J]=efu(p,u);  % calculate H 
v=u(1:p.np); w=u(p.np+1:2*p.np); % extract soln-components, states 
l1=u(2*p.np+1:3*p.np); l2=u(3*p.np+1:4*p.np); % co-states 
f1=(g*w.*v.^eta-d*(1+del*v)).*v-h; f2=rp*(beta+xi*v)-(up*v+rw).*w; % f1,f2
f3=rho*l1-pp*al*h./v-l1.*(g*(eta+1)*w.*v.^eta-2*d*del*v-d-al*h./v)-l2.*(rp*xi-up*w); 
f4=rho*l2-l1.*(g*v.^(eta+1))-l2.*(-up*v-rw); f=[f1;f2;f3;f4]; 
r=p.mat.K*u(1:p.nu)-p.mat.M*f; % the residual 
\end{Verbatim}
\vs{-6mm}\\\hline\vs{-4mm}
\begin{Verbatim}[commentchar=!, numbers=left] 
function jc=vegjcf(p,u); [e,h,jc]=efu(p,u); % J_c for vegOC, here just an interface
\end{Verbatim}
\vs{-6mm}\\\hline\vs{-4mm}
\begin{Verbatim}[commentchar=!, numbers=left] 
function [e,h,J]=efu(p,varargin) % extract [e,h,J] from p or u
if nargin>1 u=varargin{1}; else u=p.u; end 
par=p.u(p.nu+1:end); cp=par(11); pp=par(12); al=par(13); 
v=u(1:p.np); l1=u(2*p.np+1:3*p.np); 
gas=((pp-l1)*(1-al)./cp).^(1/al); e=gas.*v; h=v.^al.*e.^(1-al); 
J=pp*v.^al.*e.^(1-al)-cp*e;
\end{Verbatim} 
\vs{-6mm}\\\hline
\end{tabular}}
\vs{0mm}
\caption{{\small Implementation of the rhs $G$ and $J_c$ for \reff{vegcs}, 
and the aux function {\tt efu}. }\label{vtab2}}\end{table}

\begin{figure}\small
(a)\hspace{70mm}(b)\hspace{35mm}(c)\\[-3mm]
\ig[width=68mm]{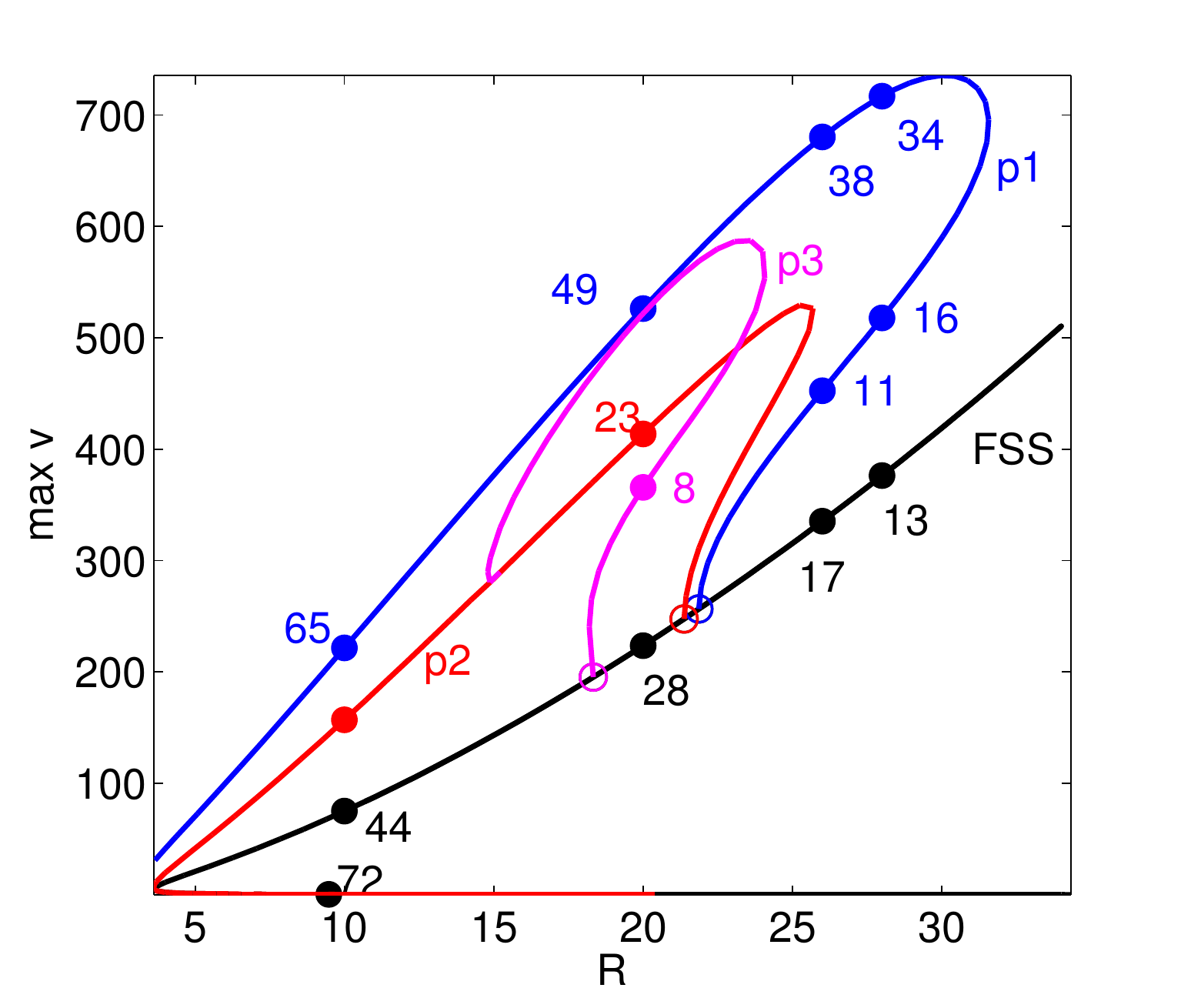}\quad\ig[width=44mm]{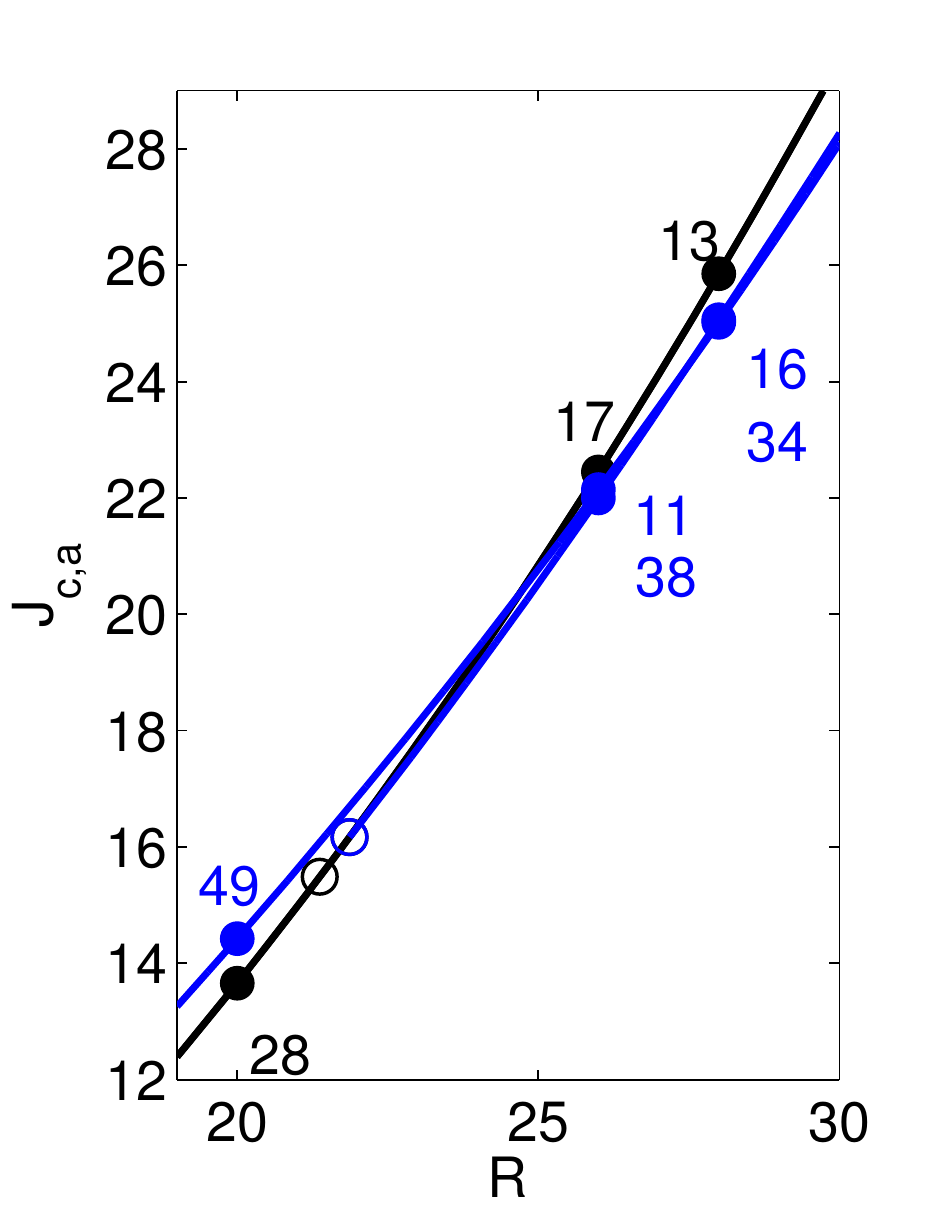} 
\raisebox{30mm}{\begin{tabular}{l}
\ig[width=35mm]{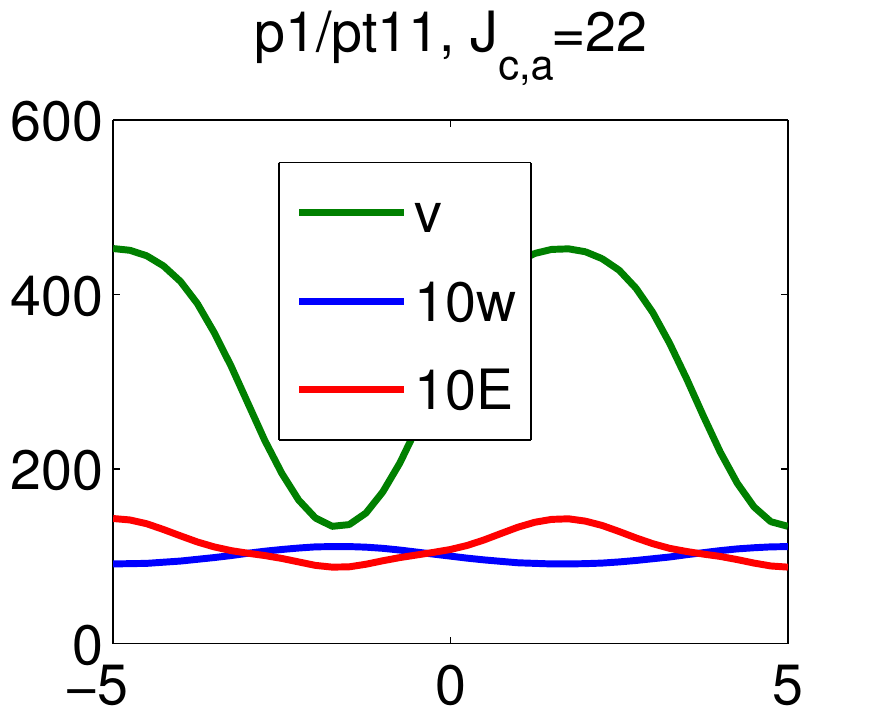}\\\ig[width=35mm]{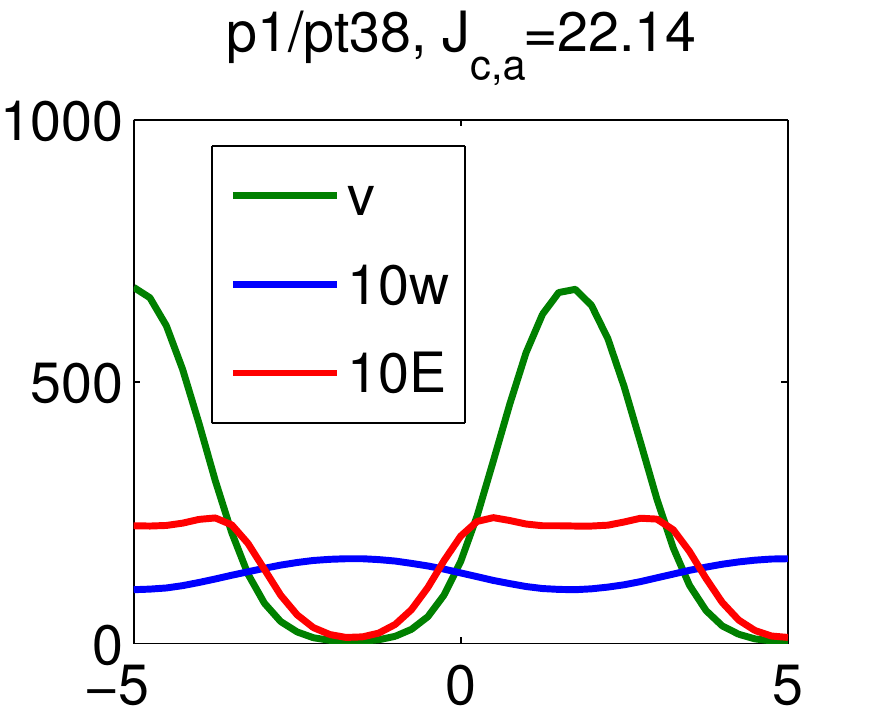}\end{tabular}}\\
\caption{{\small 
Outputs of {\tt bd1ddemo.m}. (a),(b) bifurcation diagrams of CSS 
in 1D; (c) example solutions. }
\label{vf1}}
\end{figure}

Figure \ref{vf1} shows a basic bifurcation diagram of CSS in 1D 
with $\Om=(-L,L)$, $L=5$, from the script file {\tt vegbd1d.m}, which 
follows the same principles as the one for the SLOC demo. 
 The blue branch in (a) represents 
the primary bifurcation of PCSS, which for certain $R$ have the SPP, 
and, moreover, are POSS. 
See also \cite{U15} for more plots, 
including a comparison with the uncontrolled case of so called 
``private optimization'', and 2D results for 
$\Om=(-L,L)\times (-\sqrt{3}L/2,\sqrt{3}L/2)$ yielding various POSS, 
including hexagonal patterns. 

\begin{figure}[ht]\small
(a) $R=26$, the canonical path from the lower PCSS(p1/pt11) to the FCSS 
(FSS/pt17)\\
\ig[width=41mm]{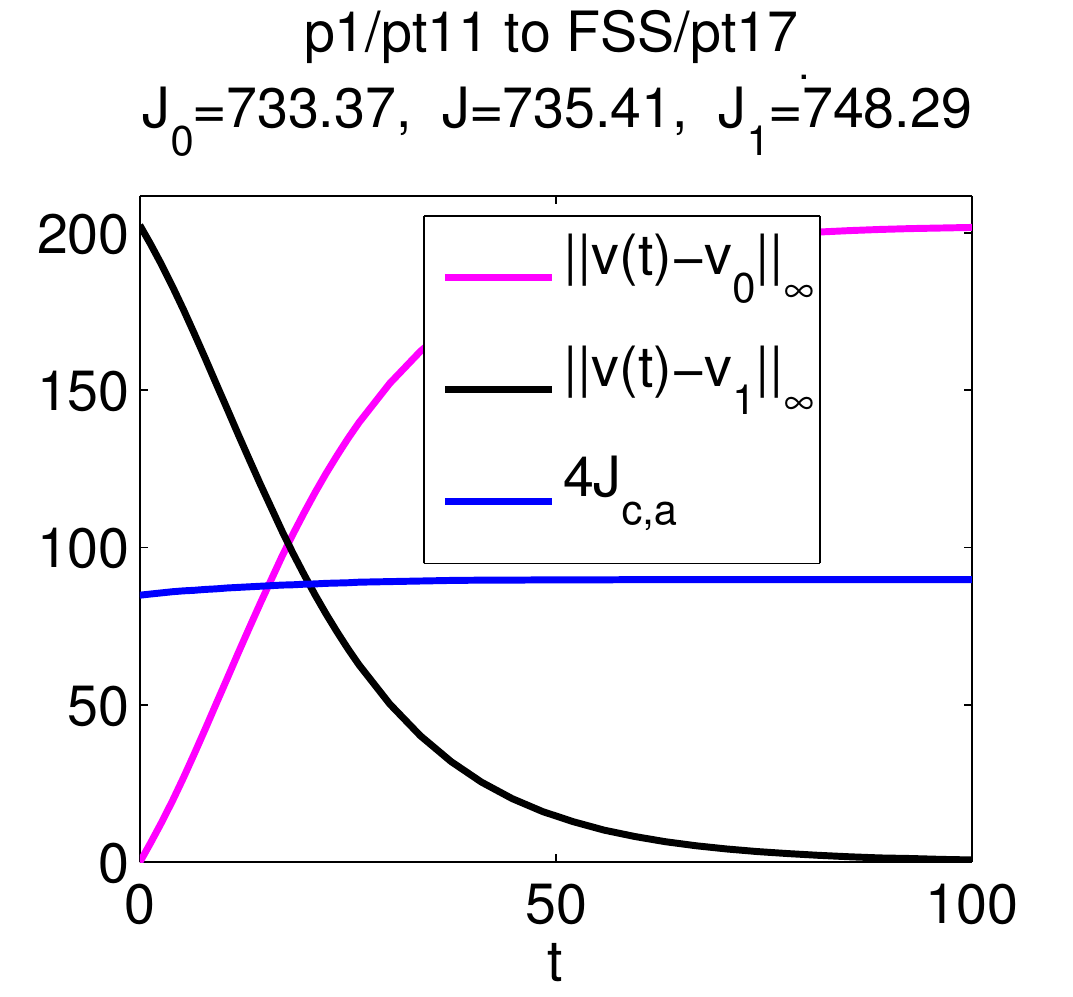}\ig[width=41mm]{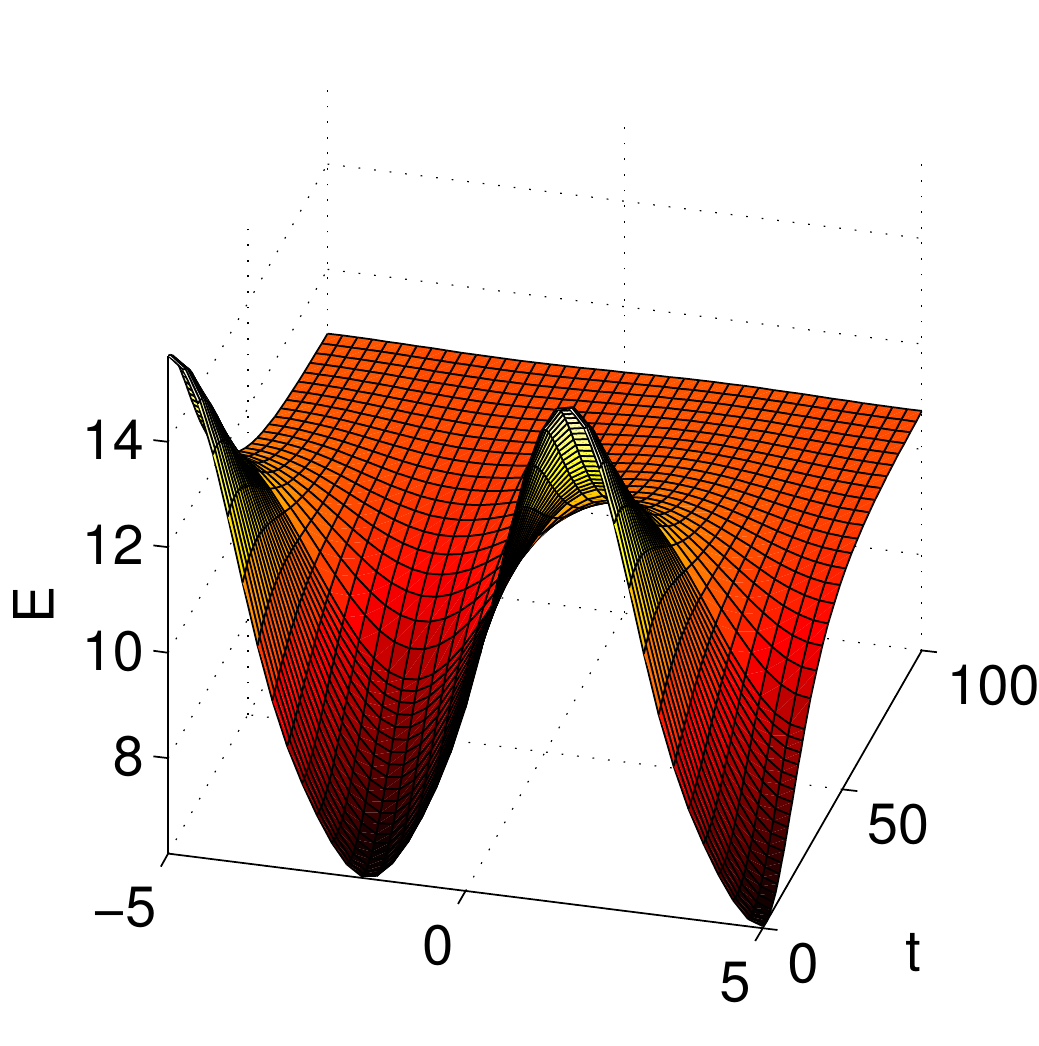}
\ig[width=41mm]{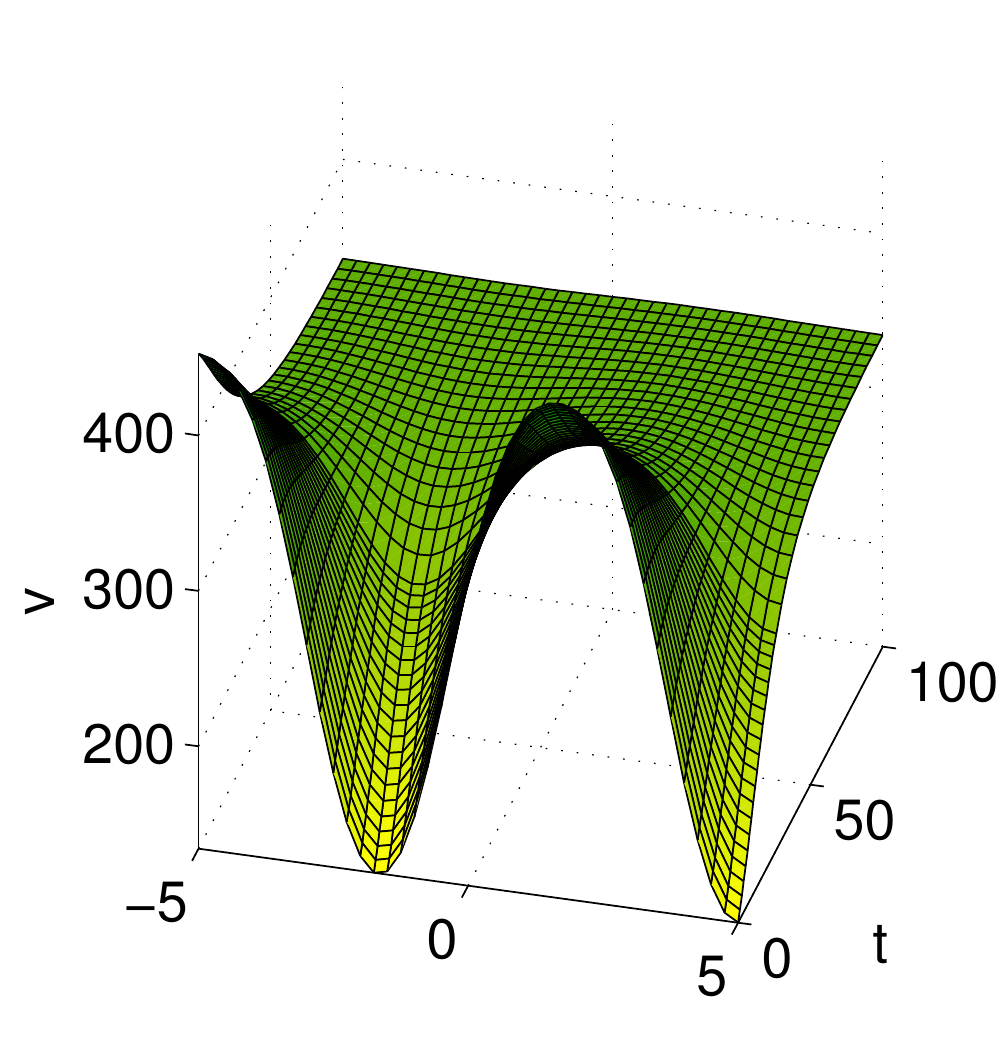}\ig[width=41mm]{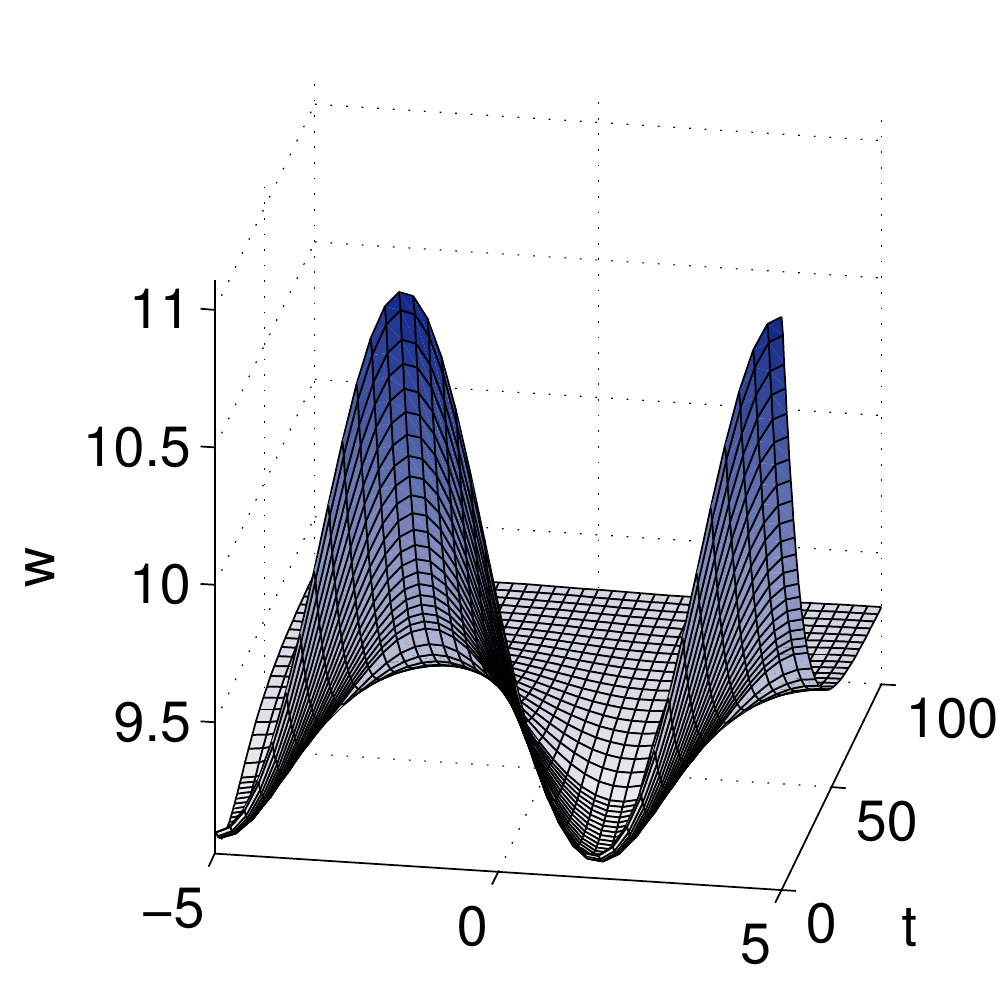}
\\[2mm]
(b) $R=26$, the canonical path from the FCSS to the upper PCSS (p1/pt38)\\
\ig[width=41mm]{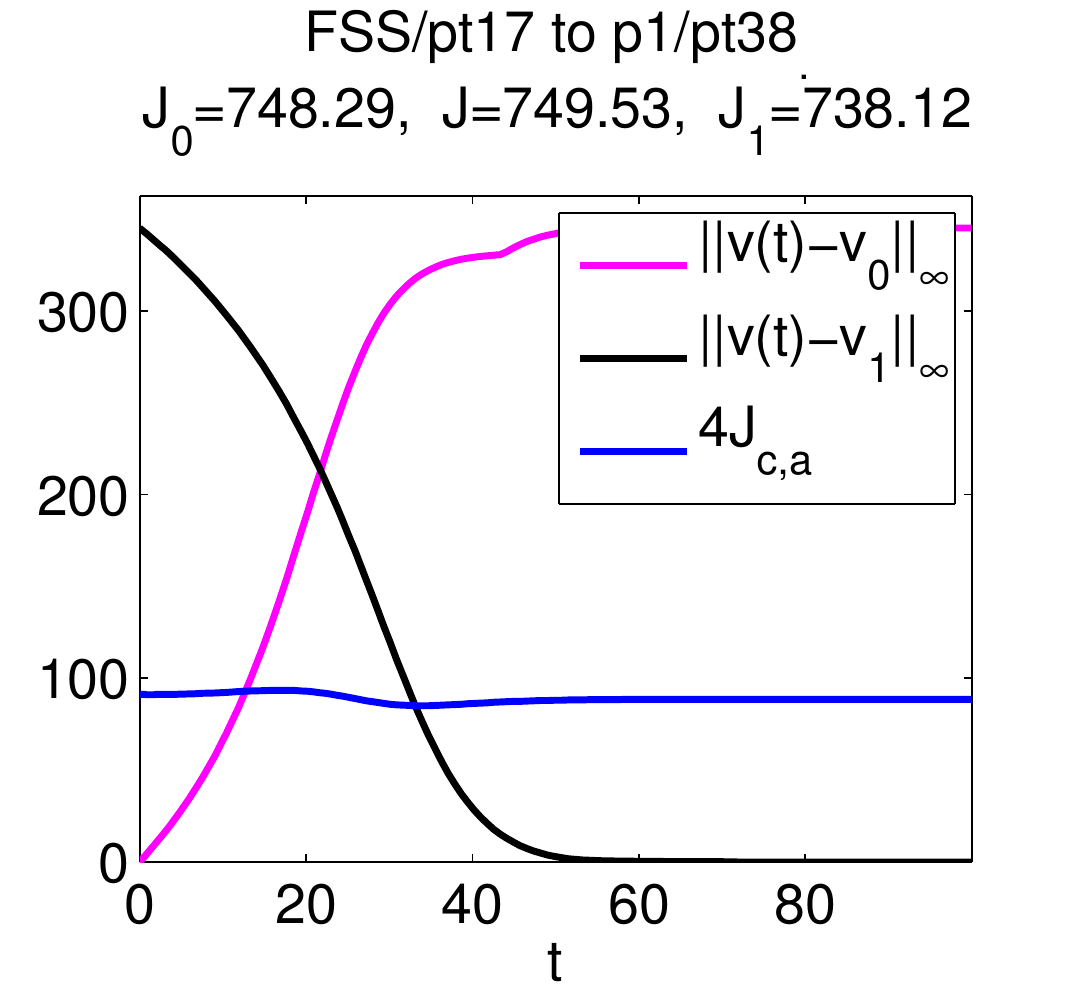}\ig[width=41mm]{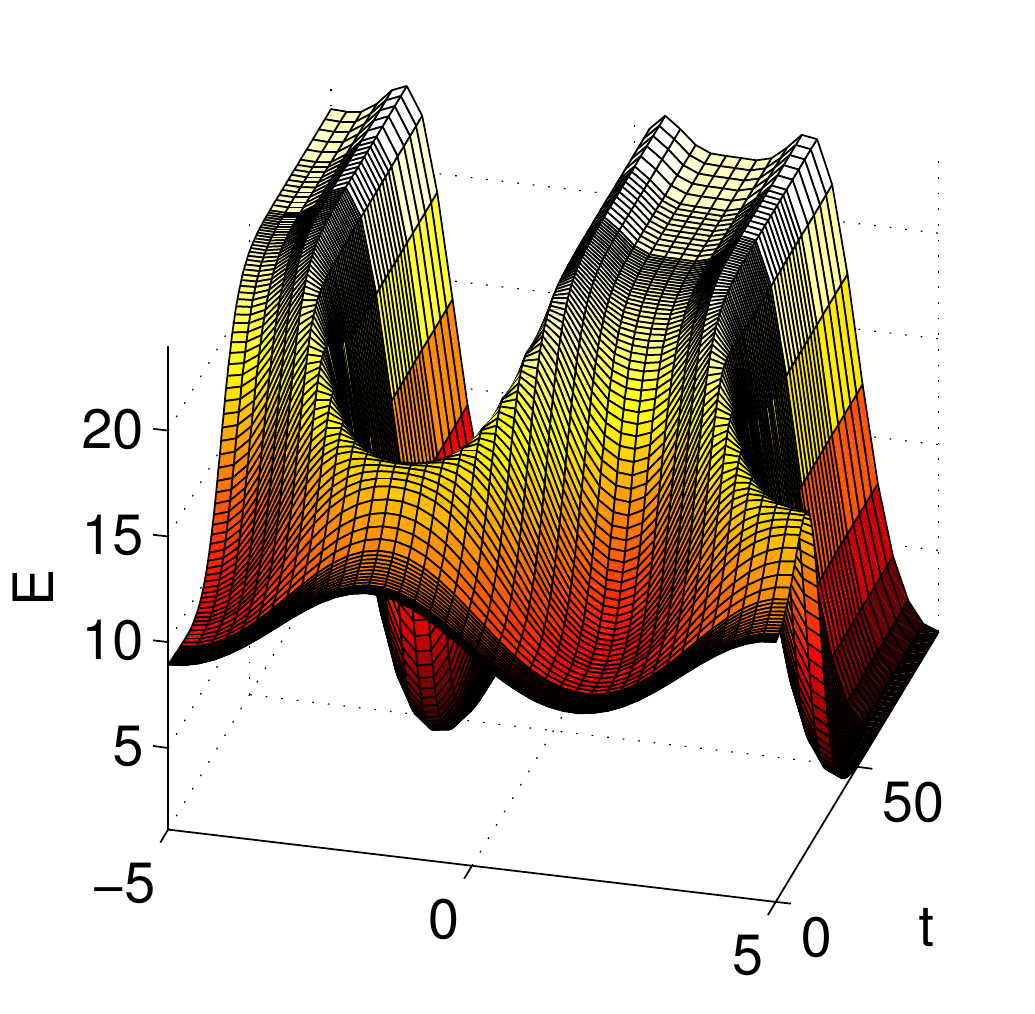}
\ig[width=41mm]{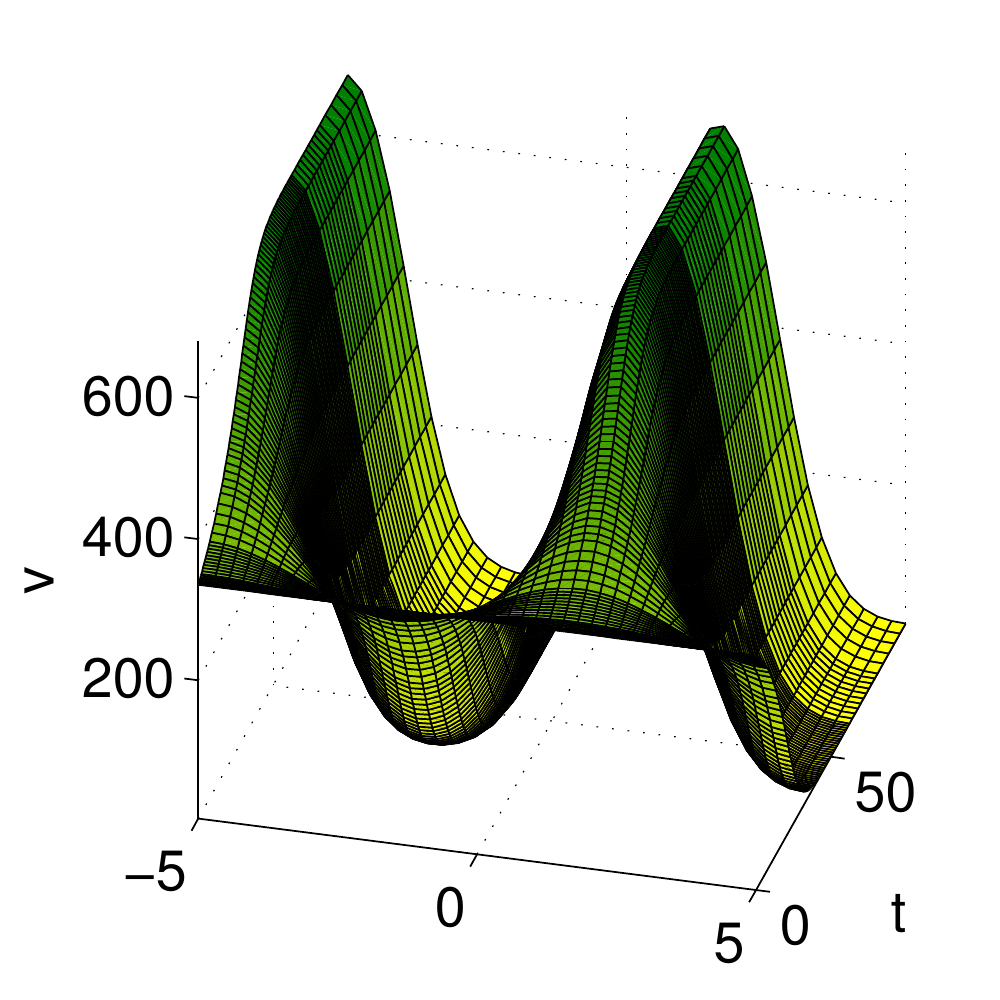}\ig[width=41mm]{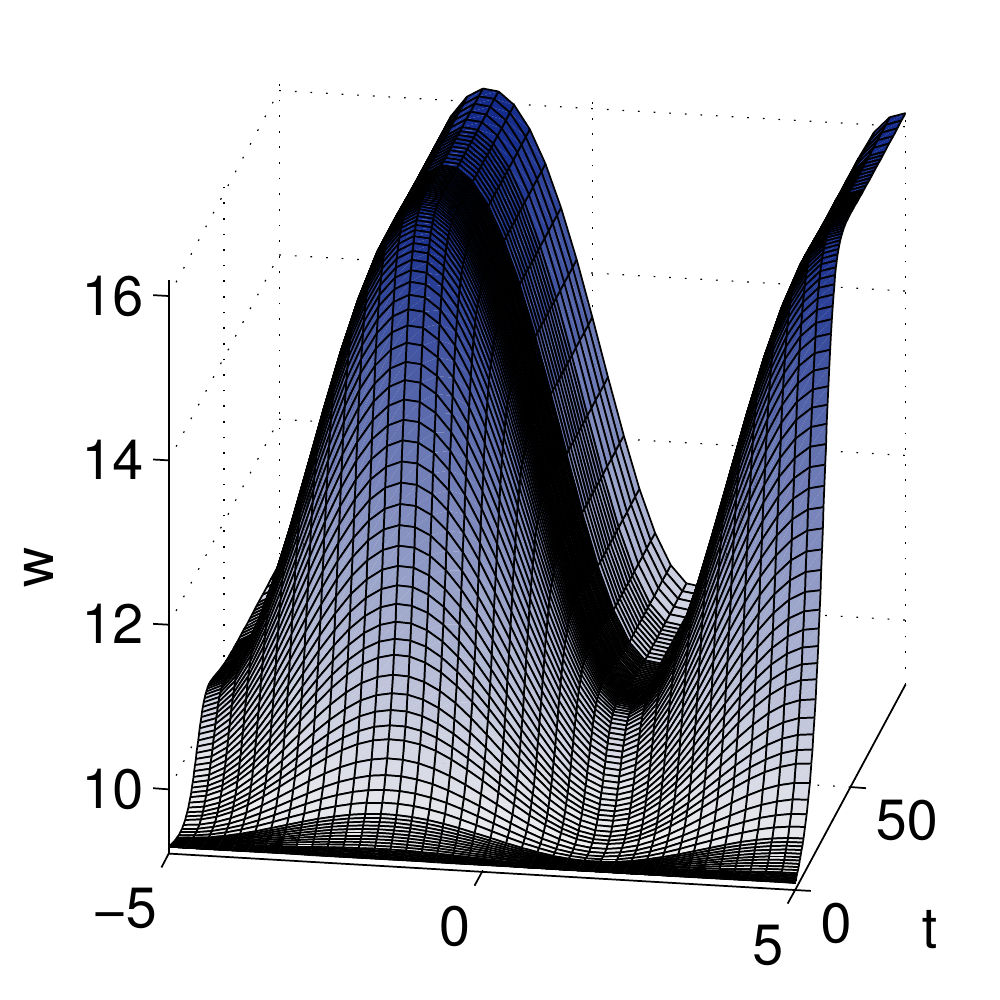}
\caption{{\small Example output of {\tt vegcpdemo.m}: Two canonical paths. 
The leftmost panels indicate the convergence behaviour, the 
current value profits, and obtained objective values. The middle and right 
panels show the strategy $E$ and the corresponding behaviour of $(v,w)$. 
See \cite{U15} for comments and more details. 
\label{vf2}}}
\end{figure}

\begin{figure}[!ht]
\bce{\small
\begin{tabular}{ll}
(a) A Skiba point at $\al\approx 0.9$&(b) Paths of (almost) equal values to the FCSS and the 
upper PCSS.\\
\ig[width=0.25\textwidth, height=35mm]{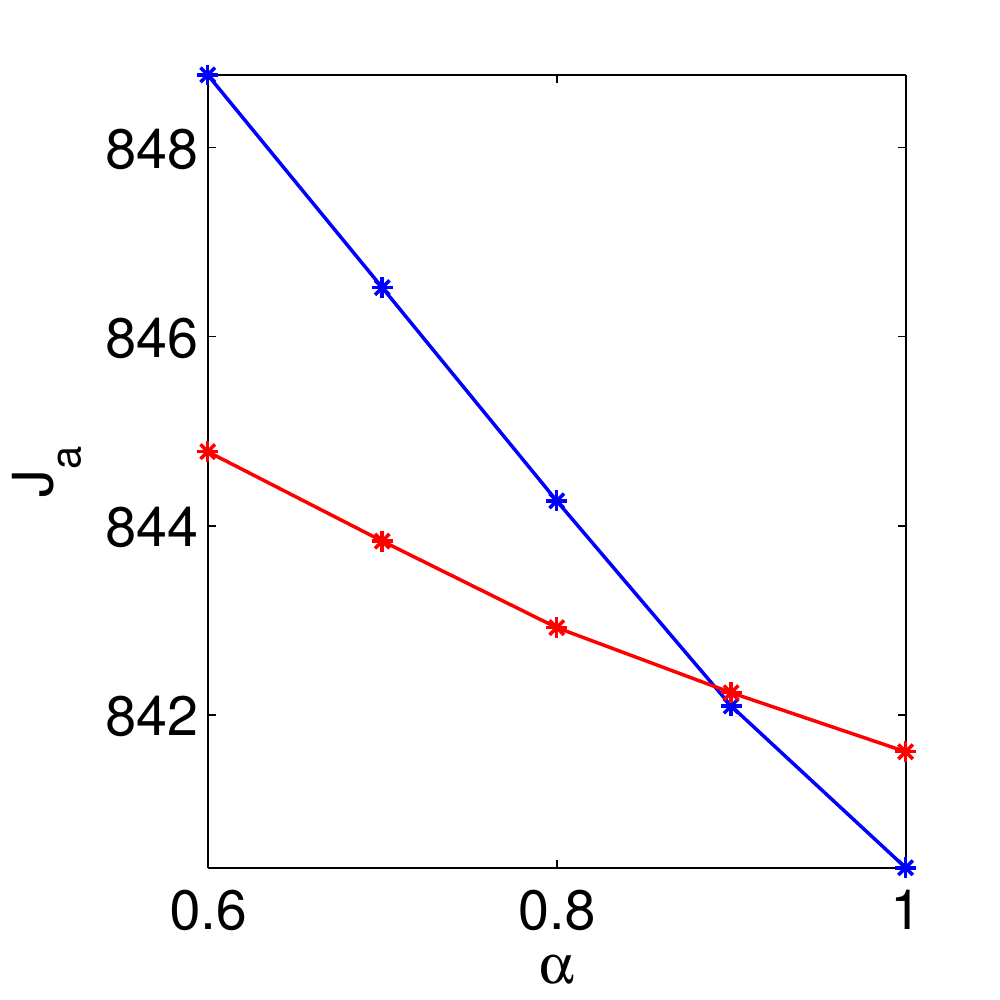}&\ig[width=0.22\textwidth]{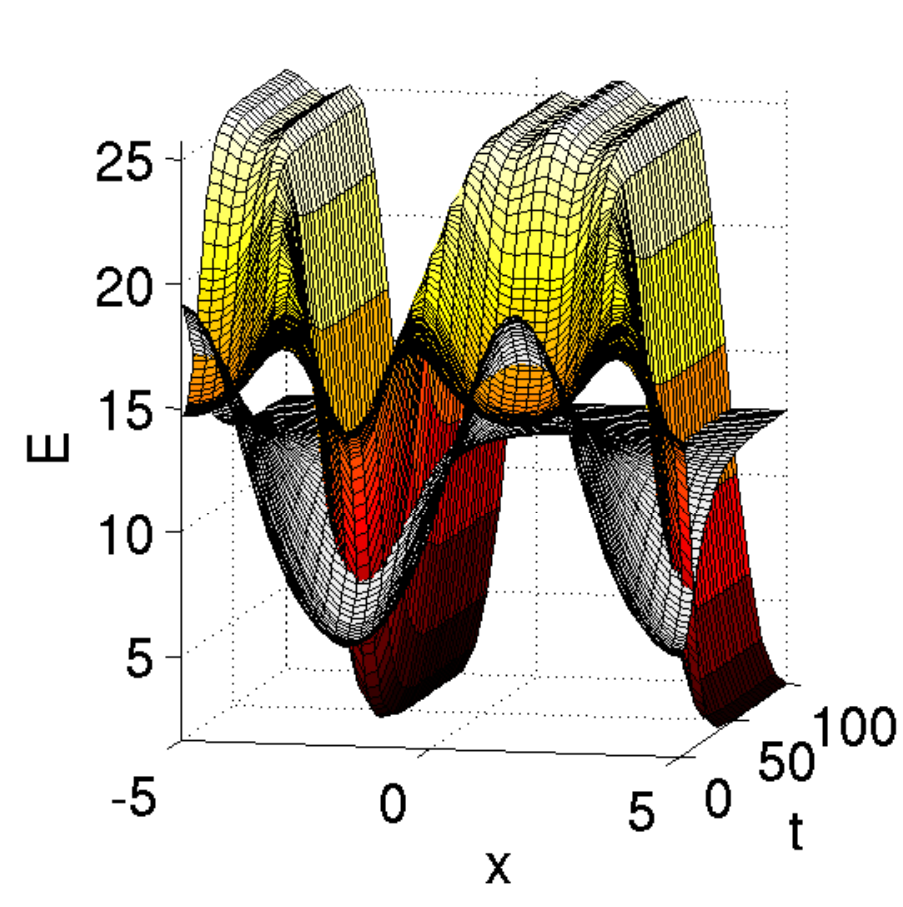}
\raisebox{0mm}{\ig[width=0.23\textwidth]{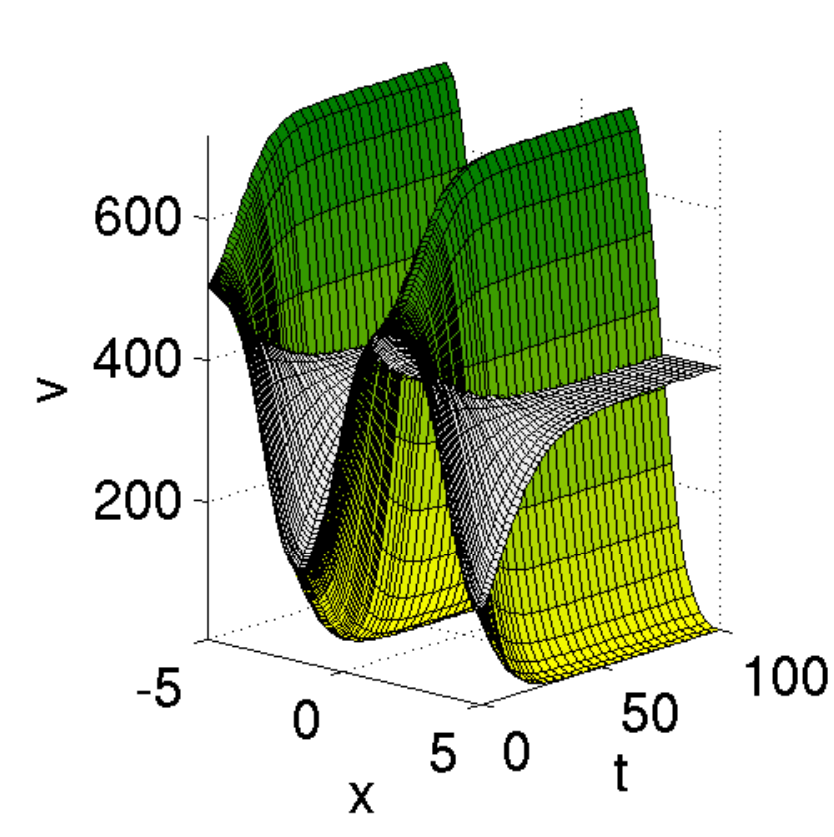}}
{\ig[width=0.23\textwidth]{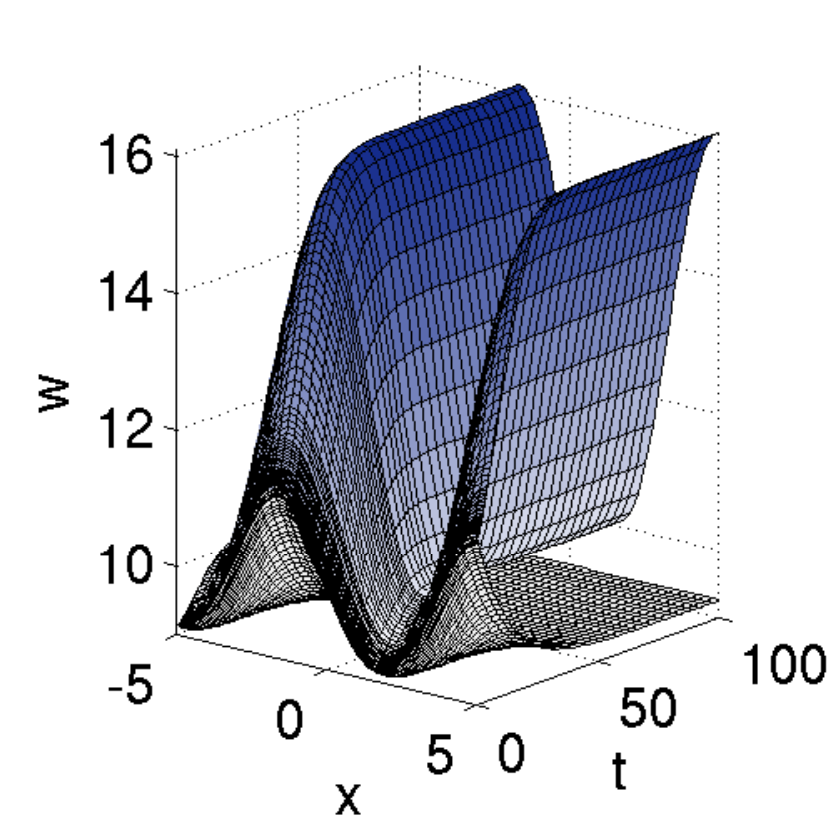}}
\end{tabular}
}
\caption{{\small $R=28$, example outputs from {\tt skibademo.m}. 
In (a), the blue 
line gives $J$ for the canonical path $t\mapsto u(t)$ 
from $(v,w)_\al(0):=\al (v,w)_{\text{PS}}+(1-\al)(v,w)_{\text{FSS}}$, 
where FSS denotes FSS/pt13, and PS denotes p1/pt16. 
The red line gives $J$ for the CP $t\mapsto \uti(t)$ 
from $P_{\al}(0)$ to the upper 
PCSS {\tt p1/pt34}. Similarly, the white surfaces in (b) 
are for $u$ and the colored ones for $\uti$. $R=28$.\label{vf3}}}
\ece
\end{figure}

The script files {\tt vegcpdemo.m} for 
canonical paths, and {\tt vegskiba.m} for a Skiba point between the 
flat optimal steady state {\tt FSS/pt13} and the POSS {\tt p1/pt34}, 
again follow the same principles as in the the SLOC demo. See 
Figures \ref{vf2}, and \ref{vf3} for example outputs, and \cite{U15} 
for a detailed discussion. In a nutshell, we find that: 
\bci
\item[{\bf (a)}] For large $R$ the FCSS is the unique CSS of \reff{vegcs}, 
and is optimal, hence a globally stable FOSS (Flat Optimal Steady State).  
\item[{\bf (b)}] For smaller $R$ there are branches of (locally stable) POSS 
(Patterned Optimal Steady States), which moreover 
dominate all other CSS. 
\item[{\bf (c)}] For the uncontrolled problem, Flat Steady States (FSS) 
only exist 
for much larger $R$ than the FCSS under control. 
\item[{\bf (d)}] At equal $R$, the profit $J$ (or equivalently the 
discounted value $J_c/\rho$) of the uncontrolled FSS is much lower 
than the value of the FCSS under control.
\eci

\section{Summary and outlook}\label{dsec}
With \poc\ we provide a toolbox to study OC problems of class \reff{oc1} 
in a simple and convenient way, in 1D and 2D. The class \reff{oc1}  
is quite general, 
and with the \pdep\ machinery we have a rather powerful tool to study the 
bifurcations of CSS. The computation of canonical paths is comparatively 
more involved. 
Essentially, our step (b) implements for the class \reff{oc1} (parts of) the 
methods explained for ODE problems in \cite[Chapter 7]{grassetal2008}, 
and implemented in 
{\tt OCMat} \url{orcos.tuwien.ac.at/research/ocmat_software/}, 
see also \cite{grass2014} for an extension of {\tt OCMat} to 
1D systems of class \reff{oc1}. In a somewhat more general sense, 
step (b) is a special case (for PDEs) of the ``connecting orbit method''. 
See \cite{auto, BPS01} and the 
references therein for earlier work on connecting orbits in ODE problems, 
including connecting orbits to periodic solutions, which 
for ODE OC problems may also be important as long-run optimal 
solutions, again cf.~\cite{grassetal2008}. 
Our setup for (b) is reasonably fast for up to 4000 degrees 
of freedom of $u$ at fixed time, e.g., 1000 spatial discretization 
points and 4 components, and 
up to 200 temporal discretization points, i.e., up to these values 
a continuation step in the calculation 
of a canonical path takes up to a few minutes on a desktop computer.

Of course, there is a rather large number of issues we do not address (yet). 
Besides periodic long-run optimal 
solutions, one of these are state or control 
inequality constraints that frequently 
occur in OC problems. 
For instance, in the SLOC model we need non-negativity of $P$ and $k$, 
and similarly of $v,w$ and $E$ in the vegOC model. In our examples 
we simply checked these a posteriori and found them to be always fulfilled, 
i.e., {\em inactive}. If such constraints become {\em active} 
the problem becomes much more complicated.  
Some extensions in this direction will be added as required 
by examples.  

Clearly, it is 
tempting to recombine steps (a) and (b) again, at least for 
specific purposes. One example would be the continuation of canonical paths 
in a parameter $\eta$. Naively, this could be done ``by hand'' by using 
a canonical path $u(\cdot,\eta)$ between $v_0$ (or $v_0(\eta)$) and $\uh( \eta)$ 
as an initial guess for a canonical path $u(\cdot,\eta+\delta)$ 
between $v_0$ (or $v_0(\eta+\del)$) and $\uh(\eta+\del)$, at the parameter value 
$\eta+\del$. This, however, does not directly allow to check for bifurcations 
of canonical paths, and, perhaps more importantly, requires the 
recalculation of $\Psi$ at each new $\uh(\eta)$. See \cite{BPS01,Pam01} 
for the ``boundary corrector method'' as an approach to avoid the 
latter, and, moreover, for continuation methods in the full t-BVP that 
for instance also allow the computation of Skiba-curves (cf.\S\ref{skibasec}) 
in 0D, cf.~also \cite[\S7.7--\S7.8]{grassetal2008}.

As currently {\tt p2pOC} is based on \pdep, 
it works most efficiently for spatial 2D problems, while 1D 
problems are treated as very narrow quasi 1D strips. Presently, 
\pdep\ is extended to efficiently treat also 1D and 3D problems, 
based on the package {\tt OOPDE} 
\url{www.mathe.tu-freiberg.de/nmo/mitarbeiter/uwe-pruefert/software}. 
Thus, {\tt p2pOC} will soon provide a genuine 1D setting as well. 

\renewcommand{\refname}{References}
\renewcommand{\arraystretch}{1.05}\renewcommand{\baselinestretch}{1}
\small
\newcommand{\etalchar}[1]{$^{#1}$}


\end{document}